\documentclass[12pt]{article}

\usepackage[english]{babel}

\usepackage[a4paper,top=2cm,bottom=2cm,left=3cm,right=3cm,marginparwidth=1.75cm]{geometry}

\usepackage{amsmath, amsthm, amssymb}
\usepackage{graphicx}
\usepackage[colorlinks=true, allcolors=blue]{hyperref}
\usepackage{verbatim}
\usepackage{mathtools, mathrsfs}
\usepackage{enumitem}
\usepackage[capitalize]{cleveref}
\usepackage[normalem]{ulem}
\usepackage{svg}
\usepackage{natbib}

\theoremstyle{definition}
\newtheorem{definition}{Definition}
\newtheorem{theorem}{Theorem}
\newtheorem{proposition}{Proposition}
\newtheorem{lemma}{Lemma}
\newtheorem{algorithm}{Algorithm}

\theoremstyle{remark}
\newtheorem{corollary}{Corollary}
\newtheorem{remark}{Remark}
\newtheorem{example}{Example}

\newtheorem{fact}{Fact}
\newtheorem{conjecture}{Conjecture}

\newcommand{\R}{\mathbb{R}}
\newcommand{\G}{\mathcal{G}}

\title{Constructing embedded surfaces for cellular embeddings of leveled spatial graphs
}
\author{S. Barthel, F. Buccoliero}

\begin{document}

\maketitle
\begin{abstract}
For a given spatial graph $\G \subset \mathbb{R}^3$, we would like to find a closed orientable surface $\mathcal{S}$ embedded in $\mathbb{R}^3$ in which $\G$ is cellular embedded. However, for general $\G$ this is not possible. We therefore define a property of spatial graphs, called leveled, to show that for leveled spatial graphs with a small number of levels, a surface $\mathcal{S}$ can always be found. The argument is based on decomposing $\G$ into spatial subgraphs that can be placed on a sphere and on cylinders attached as handles, in such a way that the resulting surface contains a cellular embedding of~$\G$.

We generalize the procedure to an algorithm that, if successful, constructs $\mathcal{S}$ for leveled spatial graphs with any number of levels. 
We conjecture that all connected leveled embeddings can be cellular embedded with the presented algorithm.

\end{abstract}

\section{Introduction}

Embeddings of graphs on surfaces are studied as central objects in a wide variety of fields. They are the main object of study in topological graph theory, where existence and properties of embeddings of given abstract graphs in abstract surfaces are considered, as well as in geometric graph theory, where the surfaces are endowed with a metric and the graphs are usually assumed to have straight-line edges. In knot theory, graph embeddings are used to produce invariants for knots, such as the Heegaard genus of a knot, which is the smallest genus of a handlebody the knot can be embedded in. In synthetic chemistry, graph embeddings are frequently used to describe the bond network of molecules.

Cellular embeddings, as a natural generalization of planar graphs, have in particular been studied in topological graph theory and combinatorics. They provide a combinatorial description of graph embeddings in surfaces~\cite{stahl1978}, and they are at the core of long-standing open conjectures such as the circular (or strong) embedding conjecture~\cite{JAEGER19851}.

It is possible to cellular embed any abstract graph in a closed orientable surface. For a given spatial graph, it is still always possible to find a closed orientable surface it embeds in by considering its regular neighborhood and pushing the spatial graph on the obtained surface. But in contrast to abstract graphs, it is not possible to find a cellular embedding for every spatial graph. For example, for nontrivial knots, there exist no such surface where the knot cellular embeds. This motivates the question for which spatial graphs a closed orientable surface embedded in $\mathbb{R}^3$ can be found, such that the graph cellular embeds in it, and how the surfaces can be constructed. 

Given that graphs on surfaces have mainly be studied in topological graph theory that, in contrast to knot theory, traditionally considers all possible embeddings of an abstract graph, the scenario involving specific realizations of graphs within three-dimensional space, where surfaces must not intersect, remains underexplored. This is despite the setting being most relevant for applications taking place in the 3-dimensional Euclidean world. For example,  building a graph-like structure using self-assembly of DNA molecules by origami folding requires first finding an unknotted A-trail through the desired spatial graph \cite{ELLISMONAGHAN201769}, \cite{MonaghanTorusGraphs}. When the graph is eulerian and cellular embedded in a surface, an unknotted A-trail can be found by following the boundary of the disks \cite{ELLISMONAGHAN201769}.

Since we do not require the surface to be standardly embedded (e.g. a torus is allowed to be knotted), it is in principle possible to obtain all cellular embeddings of spatial graphs by starting from cellular embeddings of free bouquets of circles, on which subdivisions of edges and splittings of vertices are performed on the surface. However, our goal is to construct a surface for a given spatial graph. Note that rotation systems, although equivalent to cellular embeddings for abstract graphs \cite{stahl1978}, do not ensure the described surface to be free of self-intersections.

A natural approach to construct a surface for a given spatial graph is to glue together non-intersecting properly embedded disks that are bounded by some cycles of the graph. A sufficient albeit unnecessarily strong condition for the existence of such disks is $k$-flatness of the spatial graph \cite{NikkuniOnNS}. Indeed, it suffices to find a collection of cycles that cover all edges of the graph such that each edge is covered at most twice, and that bound disks satisfying the property above. 
The union of the disks and the spatial graph produces a cellular embedding of the spatial graph in the resulting surface. In case that every edge belongs to exactly two different cycles, the embedding is a circular (or strong) embedding. However, finding a collection of such disks is in general an intractable problem, as it includes detecting the unknot, non-trivial links with any number of components, Borromean links, and more complicated entanglements such as nonfree knot- and linkless spatial graphs such as the Kinoshita theta-curve. 

 
We therefore follow another approach to construct the surface. Without loss of generality, we can restrict to connected spatial graphs. Indeed, if distinct connected components are nontrivially linked, the spatial graph cannot be cellular embedded, and unlinked components can be considered independently. The idea is to start with a subgraph of the spatial graph that cellular embeds in the sphere, and then iteratively attach handles on which the remaining parts of the graph sit on. If the graph follows a longitude and a meridian of every handle, the resulting embedding is cellular. Since in our construction a longitude of the cylinder that generates a new handle is always followed by some subgraph, cellularity is implied if the two boundary components of the cylinder are attached to different disks of the already generated surface. If the two boundary components of the cylinder are attached to the same disk, it must be ensured that the graph also runs through the meridian of the resulting handle. Furthermore, conditions need to be found that guarantee that the construction results in an embedded surface, i.e., handles should be attached in such a way that they do not intersect the already constructed surface parts. For such a procedure to work, we require the spatial graph to have some properties and we call such spatial graphs \emph{leveled}. They can be viewed as a generalization of planar graphs and are distinct from but related to open book embeddings~\cite{leveled_Properties}.

We show in \cref{prop: few levels can always be cell embedded} that all leveled spatial graphs with up to four levels can be cellular embedded in a surface that itself is embedded in $\mathbb{R}^3$ by presenting a constructive argument. We give a combinatorial description for hamiltonian leveled spatial graphs which allows to describe an algorithm whose successful termination constructs a cellular embedding of the spatial graph in \cref{algo: hamiltonian leveled}. We extend the algorithm to non-hamiltonian leveled graphs in \cref{algo: non-hamiltonian algorithm}, and finally generalize it to \emph{multi-leveled} spatial graphs, that are compositions of several leveled spatial graphs. 

We give examples where the algorithm fails to find a cellular embedding for a given starting position. However, in all such examples another initiation can be chosen such that the algorithm successfully terminates. We conjecture that this is always the case. 


\section{Preliminaries}
All graphs $G$ are undirected finite graphs. Graphs are allowed to have loops and parallel edges. The \emph{fragments} (or bridges) of a graph $G$ with respect to a cycle $C$ of $G$ are the closures of the connected components of $G - C$~\cite{MR81471}. Two fragments are \emph{conflicting} if they either have pairs of endpoints on $C$ that alternate, or have at least three endpoints in common \cite{MR1367739}. 
The \emph{conflict graph} of a graph $G$ with respect to a cycle $C$ is the graph that has a vertex for each fragment of $G$ with respect to $C$, and two vertices are adjacent if and only if the corresponding fragments conflict.
A connected graph is \emph{hamiltonian} if it contains a cycle that visits every vertex exactly once, i.e.,  all its fragments with respect to the cycle are the closures of edges or loops.
A \emph{graph
embedding} is an embedding $f : G \rightarrow S
^3$ of a graph~$G$ in $S^3$ up to ambient isotopy. The corresponding \emph{spatial
graph $\G$} is the image of this embedding. A spatial graph is \emph{trivial} if it embeds in~$S^2$. An abstract graph is \emph{planar} if it admits a trivial embedding. 
A \emph{diagram} of a spatial graph~$\G$ is the image of a projection of $\G$ onto $\mathbb{R}^2$ such that all intersections in the image are transversal, no vertices and at most two points are mapped to a crossing, and each double point is assigned as an over- or under-crossing.
We focus on a special class of spatial graphs, the leveled graphs.
\begin{definition}
    A \textbf{leveled spatial graph} $\G$ is a connected spatial graph that contains an unknotted cycle $C$ called its \textbf{spine} such that each fragment $f$ of $\G$ with respect to $C$ can be embedded in a disk $D_f$ whose boundary is identified with $C$,  with interior $\mathring{D_f}$ disjoint from $\G - f$ and disjoint from any other $D_g$, for $g$ another fragment of $\G$. 

\end{definition}

We define an ordering of the fragments of $\G$ as follows: fix a diagram $D$ of a leveled spatial graph $\G$ such that the spine $C$ of $\G$ does not have crossings on $D$. If the fragment~$f$ crosses over fragment~$g$, then $f>g$. Fragments that are not comparable by the partial order are said to be \emph{on the same level}. Note that conflicting fragments must be on different levels. Also, it is not possible that two fragments $f,g$ are such that $f>g$ and $g>f$. Indeed, if this was the case, the interior $\mathring{D_f}$ of $D_f$ could not be disjoint from $D_g$.

This definition allows us to assign the fragments of a leveled spatial graph to their level. The fragments that do not cross over any other fragments are at level~$1$, denoted $L_1$.

\begin{remark}\label{remark:cyclic_ordering_levels}
    If the leveled spatial graph $\G$ has $n$ levels $L_i, 1\leq i\leq n$, there is an action of the cyclic group $C_n = \langle x \mid x^n=1\rangle$ on the partition of levels $\mathcal{L} = \{L_1, \dots, L_n\}$. The action is given by sending $x \cdot L_i$ to $L_{i+1 (\mathrm{mod}\ n)}.$
\end{remark}
We say that a leveled spatial graph~$\G$ is in \emph{representative leveling} if for each fragment~$f$ at
level~$i > 1$, there exists a fragment~$g$ at level~$i-1$ such that $f$ and $g$ conflict. In other words, $\G$ is in representative leveling if every fragment belongs to the lowest possible level. Every leveled spatial graph can be put into representative leveling by assigning fragments to the lowest possible level. Note that the cyclic reordering of levels does not maintain representative leveling.

Every connected trivial spatial graph is a leveled graph with exactly one level. Consequently, all planar graphs admit a leveled embedding, but that is not true for all graphs. For example,  the graph~$K_5\_ K_5$ that is obtained by joining two copies of $K_5$ by an edge (\Cref{fig: multileveled embedding}) cannot be embedded as a leveled graph, since every choice of a cycle implies a non-planar fragment. Therefore, no spine can be found for $K_5\_ K_5$.

In this paper, all surfaces are considered to be closed, orientable, and embedded in $S^3$.
The spatial graph $\G$ embeds \emph{cellular} in the surface~$\mathcal{S}$ if the complement~$\mathcal{S} \setminus \G$ of the spatial graph in the surface is a union of open disks. Two pairs of spatial graphs on surfaces $(\G_1, \mathcal{S}_1), (\G_2, \mathcal{S}_2)$ are \emph{equivalent} if there is an ambient isotopy of $S^3$ between $(\G_1, \mathcal{S}_1)$ and $(\G_2, \mathcal{S}_2)$ that keeps the graphs on the surfaces throughout the ambient isotopy. Note that this notion of equivalence for graph embeddings in surfaces is stronger than the one used in topological graph theory, which only requires surface homeomorphism instead of ambient isotopy.  

Without loss of generality, all spatial graphs considered in this paper are connected, since a disconnected graph $\G$ is cellular embedded in a surface $\mathcal{S}$ if and only if each component of $\G$ is cellular embedded in a distinct component of $\mathcal{S}$ and both $\G$ and $\mathcal{S}$ have the same number of components. Indeed, if two connected components~$\G'$ and $\tilde{\G}$ of $\G$ are embedded in a connected surface $\mathcal{S}$, there exists a simple closed curve $c\in \mathcal{S}\setminus \G$ which encloses $\G'$. Since $c$ does not bound a disk in $\mathcal{S} \setminus \G$, the embedding is not cellular.

Since the subgraph of a leveled graph that consists of its spine and fragments of two levels $i$ and $j$ with $i<j$ is trivial by definition, it is always possible to embed a sphere $\mathcal{S}_0$ in $\R^3$ such that the spine $C$ of $\G$ is embedded in $\mathcal{S}_0$, all fragments of level~$i$ are embedded in one of the two closed disks $\overline{\mathcal{S}_0 \setminus C}$, and the fragments of level~$j$ on the other disk.
The disk on which fragments of level~$i$ are embedded in is called the \emph{lower hemisphere $\mathcal{S}_l$} of $\mathcal{S}_0$, and the other disk is called the \emph{upper hemisphere~$\mathcal{S}_u$}. 

\begin{lemma}\label{Lemma:3degree}
    Let $\G$ be a spatial graph with a vertex~$v$ that is the endpoint of two fragments $f_1, f_2$ with respect to a cycle $C$. It is always possible to split a vertex $v$ into two vertices $v_1$ and $v_2$ by adding an edge to $C$ such that $v_1$ is an endpoint of $f_1$ and $v_2$ is an endpoint of $f_2$, and the fragments in the new graph conflict if and only if they conflict in $\G$. Furthermore, the splitting can be chosen such that for any diagram of $\G$ the crossings remain unchanged.
\end{lemma}
\begin{proof}
A vertex $v$ of degree $k$ can always be split by replacing a small disk around $v$ with a new disk that leaves the intersection of the graph with the boundary of the disk unchanged, and replaces the vertex $v$ with two vertices $v_1$ and $v_2$ of degree $m$ and $n$ respectively that are joined by exactly one edge such that $m+n=k+2$. This procedure leaves the crossings of the diagram unchanged. In the following, every vertex splitting is performed in this manner.
\\
Two fragments that share an endpoint are conflicting by definition if they either share two more endpoints, or have three more endpoints that are alternating between the fragments. In the first case, splitting $v$ results in a sequence of alternating endpoints and therefore the fragments remain conflicting. In the second case, splitting $v$ results in four vertices that alternate between being endpoints of $f_1$ and $f_2$, therefore keeping $f_1$ and $f_2$ conflicting. \\
If the two fragments that have the vertex $v$ in common are non-conflicting, following $C$ starting from $v$ runs through all endpoints of one fragment before running through all endpoints of the other. Splitting $v$ as described above does not introduce an alternation of the endpoints, and therefore keeps the fragments non-conflicting.
\end{proof}
By repeatedly splitting vertices on the spine of a leveled spatial graph in the way described in \cref{Lemma:3degree}, their degree can be reduced to three.

\begin{figure}
    \centering
    \includegraphics[width=0.9\textwidth]{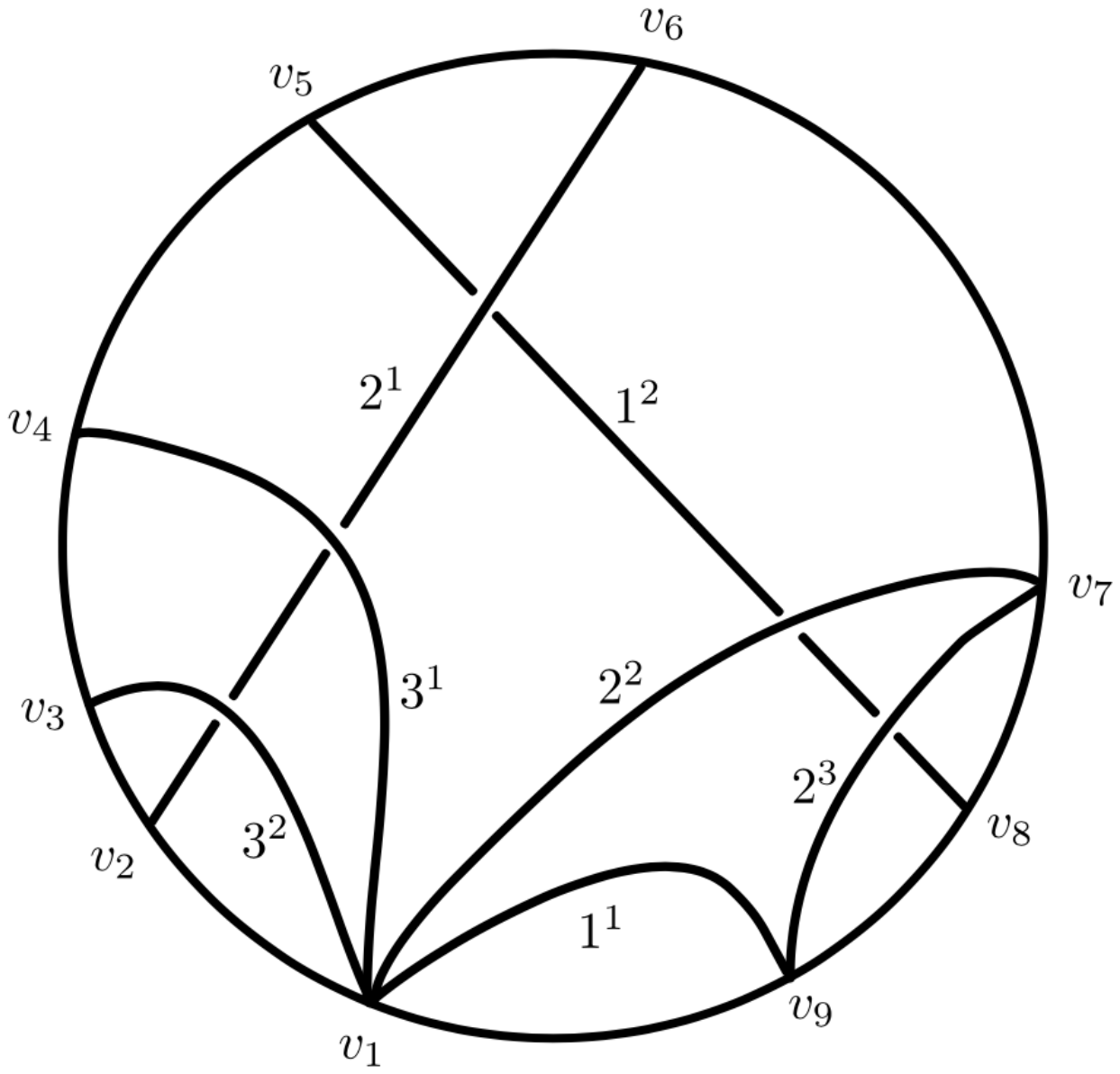}
    \caption{A hamiltonian leveled spatial graph with spine list $(1^1, 2^2, 3^1, 3^2, 2^1, 3^2, 3^1, 1^2, 2^1, 2^2, 2^3, 1^2, 2^3, 1^1)$. The superscripts distinguish fragments at same levels.}
    \label{fig:ex_string_list}
\end{figure}

A \emph{hamiltonian leveled spatial graph} $\mathcal{H}$ is a leveled spatial graph whose spine is a hamiltonian cycle of the underlying abstract graph. For constructing a cellular embedding for a hamiltonian leveled spatial graph, due to \cref{Lemma:3degree}, it is possible to restrict without loss of generality to simple graphs, i.e. graphs without loops and parallel edges. In this case, all fragments are edges with two distinct endpoints. Hamiltonian leveled spatial graphs can be completely described combinatorially by their \emph{spine list} that encodes the order in which fragments are attached to the spine together with their level. If several fragments share an endpoint, the order in the spine list is chosen such that the alternation of fragment labels is avoided. The spine list is unique up to permutation of the labels of parallel edges and cyclic permutations of its symbols, which corresponds to different choices of start points on the spine, see \Cref{fig:ex_string_list}.

For general leveled spatial graphs, spine lists descriptions can be produced analogously, introducing extra choices if the choices of splitting common vertices lead to the same minimal number of conflicts between the fragments, see \cref{Lemma:3degree}. The spine list then determines the leveled graph up to replacing fragments with other fragments that have same endpoints.

Note that the spine list of a leveled embedding $\G$ contains the information of whether or not fragments are conflicting.

Lemma~\ref{Lemma:3degree} explains how to construct the spine list of a leveled embedding that has fragments with shared endpoints. When a graph is embedded in a surface, contracting edges and splitting vertices does not change the surface. Therefore we can from now on assume without loss of generality that leveled spatial graphs have vertices of degree three on their spine, which allows the use of spine lists.

\section{Cellular embeddings of leveled spatial graphs with few levels}
In case a leveled spatial graph has a low number of levels, we give an exhaustive solution for its cellular embeddability.

\begin{proposition}\label{prop: few levels can always be cell embedded}
    If a spatial graph $\G$ has a leveled embedding with at most four levels, then there exists an oriented closed surface $\mathcal{S}$ embedded in $\R^3$ where $\G$ cellular embeds.
    \end{proposition}

    \begin{proof}
       
Assume $\G$ is in representative leveling. Let $L_i$ be the union of all fragments of level~$i$, and let $L_{1,2}:= C\cup L_1 \cup L_2$ be the subgraph of $\G$ that consists of the spine $C$ of $\G$, the fragments of level~$1$ and the fragments of level~$2$. As $L_{1,2}$ is trivial by definition, there exists an embedded sphere $\mathcal{S}_0$ where $L_{1,2}$ is embedded. Since embeddings of connected spatial graphs in a sphere are cellular, the resulting embedding $(L_{1,2}, \mathcal{S}_0)$ is cellular.

Now handles are attached to the upper hemisphere $\mathcal{S}_u$ of $\mathcal{S}_0$ to construct a surface where $L_{1,2}\cup L_3$ is embedded. For that, consider the set of disks $D_u$ that are the connected components of $\mathcal{S}_u \setminus \G$. Since each endpoint $p$ of a fragment in $L_3$ lies on the boundary $\partial D_p$ of at least one disk $ D_p \in D_u$, an assignment $a(p) = D_p \in D_u$ with $p \in \partial D_p$ can be chosen.
        
The assignment $a$ induces a clustering of $D_u$ in which two disks $D_p,D_{p'}$ belong to the same cluster $Cl_i$ if and only if there exists a fragment in $L_3$ with two endpoints $p$ and $p'$ such that $a(p)=D_p$ and $a(p')=D_{p'}$. 


Denote by $F_i$ the set of fragments of $L_3$ with endpoints in $Cl_i$.
Since $F_i$ is trivial, there is a sphere $S_{F_i}$ where $F_i$ embeds. For each set of endpoints $p_k$ of $F_i$ for which $a(p_k) =D_{p_k}$, remove a disk $P_{p_k}$ from $S_{F_i}$ such that $p_k \in \partial P_{p_k}$ for all $k$ and $q \notin P_{p_k}$ for any endpoint $q \neq p_k$. This construction yields an $n_i$-punctured sphere $\mathcal{S}_i$ where all fragments of $Cl_i$ are embedded.

The embedding of $F_i$ in $\mathcal{S}_i$ is cellular, because the spatial graph $F_i$ together with the boundary cycles of the disks $P_{p_k}$ is a connected trivial spatial graph.

Since the punctured spheres $\mathcal{S}_i$ are constructed to follow the non-conflicting fragments of level~$3$, it is implied that the $\mathcal{S}_i$ are not intersecting. Removing the interiors of the disks that are in the image of $a$ from $\mathcal{S}_u$ allows us to construct an oriented closed non-intersecting surface. That is done by identifying their boundaries with the boundaries of the punctures of $\mathcal{S}_i$ that respect the induced orientations on the boundaries of the surfaces.  


Analogously, $L_4$ is (cellular) embedded in $n_j$-punctured punctured spheres $\mathcal{S}_j$ that are attached to the lower hemisphere $\mathcal{S}_l$. 

The punctured spheres where $L_4$ embeds do not intersect the surface where $L_{1,2}\cup L_3$ is embedded. That is due to there being four levels only: each level is attached to only one of the hemispheres and the crossing types between the levels are of only one type, making the attachment of handles to the upper hemisphere and to the lower hemisphere independent of one another. The result is an embedding of $\G = L_1 \cup L_2 \cup L_3\cup L_4$ in a surface of genus $\sum_{i} (n_i -1) + \sum_{j} (n_j -1) $.

It remains to be seen that the embedding is cellular. The embedding of $\G$ in the resulting surface is cellular since each $\mathcal{S}_k$ is attached to different disks of $\mathcal{S}_0 \setminus \G$, and $F_k$ is cellular embedded in $\mathcal{S}_k$, and $L_{1,2}$ is cellular embedded in $\mathcal{S}_0$.

    \end{proof}

\begin{corollary}
    The genus of the surface $\mathcal{S}$ obtained in \cref{prop: few levels can always be cell embedded} depends on the leveled embedding itself and the definition of the function $a(p)$ in the proof. There is a unique $a$ if fragments have endpoints of degree three only. For degrees greater than three, $a$ is uniquely defined if and only if neither fragments of level 3 and fragments of level 2 nor fragments of level 4 and fragments of level 1 have endpoints in common.
\end{corollary}

The restriction to four levels in \cref{prop: few levels can always be cell embedded} assures both the absence of self-intersections of the resulting surface and the cellularity of the resulting embedding. The construction of the proof of \cref{prop: few levels can always be cell embedded} ensures neither property for a leveled embedding with five or more levels, as can be seen in the following example.

\begin{example}\label{exm: graph not embeddable with the algorithm unless you change spine}
    \quad Consider the leveled embedding $\G$ shown left in \Cref{fig: changing the spine to hamiltonian it can still be cellular embedded}, that is given by the spine list $(1^1, 1^2, 1^3, 2^1, 2^2, 2^3, 3^1, 3^2, 3^3, $ $4^1, 4^2, 4^3, 5^1, 5^2, 5^3, 1^3, 1^2, 1^1, 2^3, 2^2, 2^1,$ $3^3, 3^2, 3^1, 4^3, 4^2, 4^1, 5^3, 5^2, 5^1).$
    We can apply the algorithm in the proof of \cref{prop: few levels can always be cell embedded} to embed any four levels in a surface $\mathcal{S}$; without loss of generality choose the levels~1, 2, 4, 5. This leaves level~3 to be embedded. The endpoints of fragments of level~3 lie on the boundary of two disks of $\mathcal{S}\setminus \{L_1\cup L_2 \cup L_4 \cup L_5\}$. Denote by $D_u$ (respectively $D_l$) the disk which has interior points in common with the upper (resp. lower) hemisphere $S_u$ (resp. $S_l$) of the sphere that is produced in the first step of constructing $\mathcal{S}$. Extending $\mathcal{S}$ to accommodate $L_3$ cellularly requires a single cylinder $P$ to be placed under fragments $3^1, 3^2, 3^3$ such that $P$ is attached to $\mathcal{S}$ with one end glued to $D_u$ and the other to $D_l$. 
    However, the resulting surface is non-orientable. This can be seen as follows: the order in which the three fragments are attached to the spine defines an orientation of the spine. As the three fragments do not conflict in $\G$ with respect to the spine, and $P$ is glued to $D_u$ on one end and to $D_l$ on the other, the two ends of $P$ are oriented with the same orientation as the spine. That is, if the orientation on one boundary component of the  cylinder is clockwise, it is counterclockwise for the other. The boundary orientation does therefore not correspond to the orientation that is obtained when considering the two ends of the cylinder as naturally identified. This implies that the surface obtained by gluing $P$ to the oriented surface $\mathcal{S}$ is non-orientable.\\
     This shows that the surface with $\G$ embedded is either cellular or free of self-intersections, but cannot be both.

\end{example}

\section{Constructing the surface}
\subsection{Placing cylinders under fragments}
In this section, given a spatial graph $\G$ embedded in a surface $\mathcal{S}$ formed by gluing cylinders consecutively to a sphere, we derive restrictions on the placement of fragments on the cylinders forming $\mathcal{S}$ such that the embedding of $\G$ in $\mathcal{S}$ is leveled.

Unless specified otherwise, we suppose that the spine of a leveled embedding is a hamiltonian cycle of the underlying graph. Under this assumption, every fragment is an edge of the graph. To denote the spine list of leveled spatial graphs, we use letters as variables for fragments.

\begin{figure}
    \centering
    \includegraphics[width=0.8\linewidth]{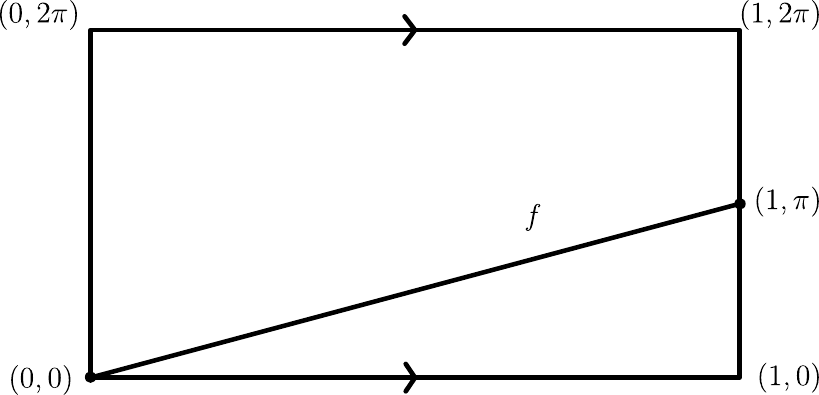}
    \caption{A fragment $f$ with a right halftwist with respect to the cylinder $P$.}
    \label{fig: right halftwist}
\end{figure}

\begin{definition}\label{def: consecutively embedded}
    Let $\G$ be a leveled embedding in representative leveling with $n$ levels. We say that $\G$ is \textbf{consecutively embedded} if it is embedded in a surface $\mathcal{S} \cong S^2 \cup_{\varphi_1} P_1 \cup_{\varphi_2} \dots \cup_{\varphi_m} P_m$ that is built successively by gluing the ends of a cylinder~$P_i$ to $\mathcal{S}_{i-1} = S^2 \cup_{\varphi_1} P_1 \cup_{\varphi_2} \dots \cup_{\varphi_{i-1}} P_{i-1}$ by a map $\varphi_i$ such that: \begin{itemize}
        \item the spine of $\G$ is embedded as equator of $\mathcal{S}_0=S^2$ and levels 1 and 2 are embedded respectively in the lower and upper hemisphere of $S^2$;
        \item each cylinder $P_i$ carries at least one fragment of $\G$;
        \item the subgraph of $\G$ embedded in $\mathcal{S}_i$ is called $\G_i$;
        \item $\varphi_i$ glues $P_i$ to $\mathcal{S}_{i-1}$ such that the two ends of $P_i$ are glued to two subsets of one or two disks of $\mathcal{S}_{i-1} \setminus \G_{i-1}$ which are each homeomorphic to a disk and incident to the spine of $\G$.
    \end{itemize}
\end{definition}
\begin{definition}
    Let $\G$ be a leveled embedding in representative leveling with $n$ levels, consecutively embedded in a surface $\mathcal{S}_i$. Let $P_i$ be a cylinder which is glued to two not necessarily distinct disks $D_1, D_2$, such that $\mathcal{S}_i = \mathcal{S}_{i-1} \cup_{\varphi_i} P_i$. A subset of the spine list of $\G$ consisting of the endpoints of the fragments embedded in $P_i$ which lie on one end of $P_i$ is called a \textbf{spine sublist} of $P_i$.
    Reversing the order of the endpoints in one of the spine sublist results in the \textbf{disk sublists} of $P_i$. That is, to obtain the disk sublists, the boundary components of the cylinder are considered as naturally identified.
    
    \end{definition}
    \begin{definition}
    Let $f$ be a fragment embedded in a cylinder $P$ such that the endpoints of $f$ lie on both ends of $P$. Consider the rectangle representation of $P$ given by $[0,1]\times [0, 2\pi]/\sim$, where the point $(t,0)$ is identified with the point $(t, 2\pi)$ for $t\in [0,1].$  We say that $f$ has a right (resp. left) \textbf{halftwist} with respect to $P$ if $f$ is ambient isotopic with fixed endpoints to the line joining points $(0,0)$ and $( 1, \pi)$ (resp. the line joining points $(0, 2\pi)$ and $(1, \pi)$). (\Cref{fig: right halftwist})
\end{definition}

\begin{remark}
    The embedding of fragments $f_1, \dots, f_n$ of a leveled spatial graph $\G$ in a cylinder $P$ where every fragment has endpoints lying on the boundary of $P$ is cellular in $P$ if and only if at least one fragment $f_i$ has endpoints on both ends of $P$.

    Indeed, $f_i$ cuts $P$ along the longitude, therefore creating a disk. Any other fragment $f_j$ then cuts a disk into two disks. Therefore, the resulting embedding~$(f_1, \dots, f_n, P)$  is cellular. On the other hand, if the embedding~$(f_1, \dots, f_n, P)$ is cellular, every simple closed curve of $P$ parallel to the boundary of $P$ (a meridian of $P$) must intersect at least one fragment $f_i$. Therefore, there exists a fragment $f_i$ with endpoints on both ends of $P$.
\end{remark}

\begin{lemma}\label{lem: two fragments on a cylinder have at most 1 halftwist}
    Let $\G$ be a leveled spatial graph consecutively embedded in a surface $\mathcal{S}$. Let $P_i$ be a cylinder which is glued to two (not necessarily distinct) disks $D_1, D_2$ of $\mathcal{S}_{i-1} \setminus \G_{i-1}$.
    If a fragment $f$ embedded on $P_i$ has two or more halftwists with respect to $P_i$, then no other fragment with endpoints on both ends of $P_i$ has more than one halftwist with respect to $P_i$. Moreover, if $f$ has three or more halftwists with respect to $P_i$, then $f$ is the only fragment with endpoints on both ends of $P_i$.
\end{lemma}
\begin{proof}
    If there exists a pair of fragments with endpoints on both ends of $P_i$ with two or more halftwists with respect to $P_i$, then the two fragments need to cross over and under each other. Therefore, the embedding would not be leveled.
    Moreover, if $f$ has three or more halftwists with respect to $P_i$, then any other fragment with endpoints on both ends of $P_i$ needs to be parallel to it on $P_i$. Therefore, the difference of halftwists with respect to $P_i$ of the two fragments is at most 1. That is impossible due to the previous argument. 
\end{proof}

\begin{proposition}\label{prop: fragments on cylinder}
    Let $\G$ be a leveled spatial graph consecutively embedded in a surface $\mathcal{S}$. Let $P_i$ be a cylinder which is glued to two (not necessarily distinct) disks $D_1, D_2$ of $\mathcal{S}_{i-1} \setminus \G_{i-1}$. Then the disk sublists of $P_i$ must be one of the following: 
    \begin{enumerate}
    \item $(a_1, \dots, a_m, b_1, \dots, b_p, c, a_m, \dots, a_1/ d_1, \dots, d_q, c, b_1, \dots, b_p, d_q, \dots,d_1);$
    \item $(a_1, \dots, a_m, c, b_1, \dots, b_p, a_m, \dots, a_1/ d_1, \dots, d_q, b_1, \dots, b_p, c, d_q, \dots,d_1);$
    \item $(s_1, \dots, s_n, \alpha_1, \dots \alpha_k, \beta_1, \dots, \beta_h, s_n, \dots, s_1 / t_1, \dots, t_q, \beta_1, \dots, \beta_h,$ \\$ \alpha_1, \dots, \alpha_k, t_q, \dots, t_1),$
    \item $(a_1, \dots, a_m, \alpha_1, \dots, \alpha_k, b_1, \dots, b_p, \beta_1, \dots, \beta_h, a_m, \dots, a_1 / \beta_1, \dots, \beta_h, $ \\$\alpha_1, \dots, \alpha_k, b_1, \dots, b_p);$
    \item $(a_1, \dots, a_m, c, \alpha_1, \dots, \alpha_k, b_1, \dots, b_p, a_m, \dots, a_1 / \alpha_1, \dots, \alpha_k, b_1, \dots, b_p, c);$
    \item $(a_1, \dots, a_m, b_1, \dots, b_p, \beta_1, \dots, \beta_h, c, a_m, \dots, a_1 / \beta_1, \dots, \beta_h, c, b_1, \dots, b_p);$
\end{enumerate}
where the slash ``/" separates the two disk sublists of $P_i$, the letters used are variables, $c$ denotes a fragment with at least two halftwists, Greek letters denote fragments with exactly one halftwist, same Latin letters indicate fragments at same level, some letters may be absent except there must be at least one fragment having endpoints on both ends of $P_i$, sequences of Latin letters $x_1, \dots, x_r, x_r, \dots, x_1$ can be added in any point of the disk sublists.
\end{proposition}

\begin{figure}
    \centering
    \includegraphics[width=0.7\textwidth]{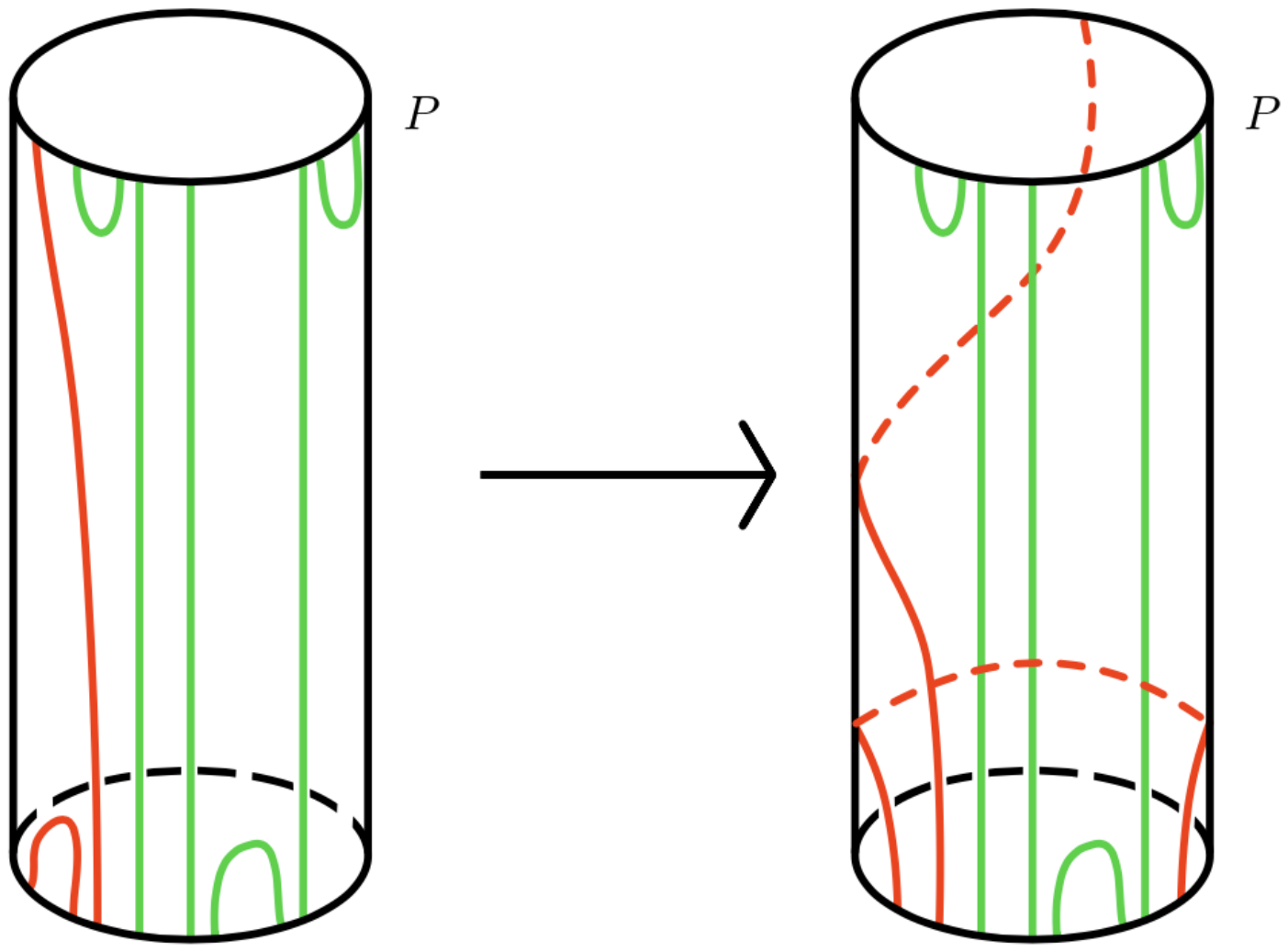}
    \caption{The endpoints of the fragments can be moved along the connected components of the boundary of the cylinder $P$ which is not attached to the surface.}
    \label{fig: shifting endpoints on cylinder}
\end{figure}

\begin{proof}
All endpoints of fragments lie on the boundary of $P_i$ by construction. The placement of fragments on $P_i$ is restricted by how fragments can be placed on a cylinder in general, and secondly by the requirement of leveledness.

For the first part, consider a cylinder $P$ that is not attached to a surface. There are two ways to embed edges in $P$ such that their endpoints lie on the boundary of the cylinder: either the two endpoints of an edge are attached to both ends of the cylinder, or both endpoints are attached to the same end. 
The endpoints of the fragments embedded in $P_i$ can be moved along the connected component of the boundary of $P_i$ to which the endpoints belong to while maintaining their embedding in $P_i$. Lemma \ref{lem: two fragments on a cylinder have at most 1 halftwist} gives conditions for the embedding to be leveled: in order to obtain the lists presented, consider the cylinder $P_i$ embedded in $\mathbb{R}^3$ but not glued to $\mathcal{S}_i$ depicted in \Cref{fig: shifting endpoints on cylinder}. Shifting the endpoints of the fragments along the boundary components of $P_i$ while respecting the conditions given by \cref{lem: two fragments on a cylinder have at most 1 halftwist} yields the six presented lists.
\end{proof}

\begin{remark}
    Lemma \ref{lem: two fragments on a cylinder have at most 1 halftwist} and \cref{prop: fragments on cylinder} restrict how fragments are embedded in a cylinder $P_i$ and the possible disk sublists of $P_i$ in a consecutively embedded leveled spatial graph.
\end{remark}

\begin{example}
    Consider the cylinder $P$ on the right of \Cref{fig: shifting endpoints on cylinder}. If we name $a_1$ the red fragment with endpoints attached to one end of $P$, $\alpha$ the red fragment with endpoints attached to both ends of $P$, $b_1, b_2, b_3$ the green fragments with endpoints on both ends of $P$ and $x_1, x_2, x_3$ the green fragments with endpoints on one end of $P$, then the disk sublists of $P$ would be $$(a_1, \alpha, b_1, b_2, x_1, x_1, b_3, a_1 / \alpha, x_2, x_2, b_1, b_2, b_3, x_3, x_3),$$ which correspond to the sublists of type $4$ of \cref{prop: fragments on cylinder}.
\end{example}

Proposition \ref{prop: fragments on cylinder} allows us to extend \cref{prop: few levels can always be cell embedded}:
\begin{corollary}\label{cor: up to 5 levels may be cell embedded}
    If a spatial graph $\G$ admits a leveled embedding with at most five levels such that at least one level consists of at most two edges, then there exists a surface $\mathcal{S}$ where $\G$ cellular embeds.
\end{corollary}
\begin{proof}
    Let the level which consists of at most two edges be level 3 of the embedding. Perform the algorithm shown in the proof of \cref{prop: few levels can always be cell embedded} on the other four levels of $\G$. This creates a surface where four levels of $\G$ are cellular embedded and the embedding has at least two disks which intersect $\mathcal{S}_0$ in different hemispheres. 
    If level 3 consists of one fragment or of two fragments whose endpoints lie in the boundaries of one pair of disks, then embed the fragments of level 3 in one cylinder; this is possible due to \cref{prop: fragments on cylinder} and the fact that the two fragments do not conflict.
    Otherwise level~3 consists of two fragment whose endpoints do not lie in the boundaries of one pair of disks. Then each fragment can be placed on a separate cylinder.

\end{proof}

\begin{figure}[h]
        \centering
        \includegraphics[width=0.5\textwidth]{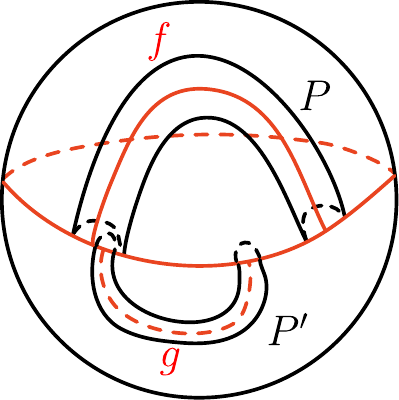}
        \caption{An upper cylinder $P$ with a fragment $f$ with upper endpoints and a lower cylinder $P'$ with a fragment $g$ with upper endpoints.}
        \label{fig: definition upper and lower endpoints and cylinders}
    \end{figure}

\begin{definition}
    Let $\G$ be a leveled embedding which is consecutively embedded in $\mathcal{S}$. Let $P$ be a cylinder in $\mathcal{S}$. We say that $P$ is an \textbf{upper} cylinder if the fragments embedded in it are at higher level than the fragments embedded in $\mathcal{S}_0$. Otherwise, $P$ is called a \textbf{lower} cylinder.
  
    The neighbourhood of the spine of $\G$ in the surface $\mathcal{S}$ has two components. We call them $N_u$ and $N_l$. For an endpoint $p$ of a fragment $f$ embedded in a cylinder $P$, let $N_p$ be a small neighbourhood of $p$ in $P$. Then $N_p$ is contained in exactly one of $N_u$ or $N_l$. We say that $p$ is an \textbf{upper endpoint} (resp. lower endpoint) if $N_p$ is contained in $N_u$ (resp. $N_l$). (\Cref{fig: definition upper and lower endpoints and cylinders})

\end{definition}

\begin{figure}[h]
    \centering
    \includegraphics[width=0.5\textwidth]{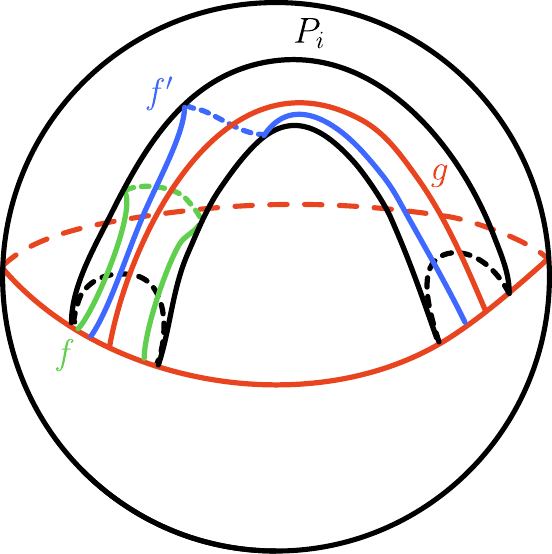}
    \caption{The situation described in cases 1 and 2 of \cref{lem: levels of fragments on cylinder}.}
    \label{fig: levels of fragments on cylinder}
\end{figure}

\begin{lemma}\label{lem: levels of fragments on cylinder}
    Consider a cylinder $P_i$ of the consecutively embedded spatial graph $\G$ in $\mathcal{S}.$
    \begin{enumerate}
        \item If $P_i$ is an upper cylinder, the fragments of $P_i$ with endpoints on only one end of $P_i$ and with at least one upper endpoint must be at lower level than the fragments they conflict with that have endpoints on both ends of $P_i$ and at least one upper endpoint in the spine sublist that contains the conflict.
        \item If $P_i$ is an upper cylinder and has a fragment with upper endpoints on both ends of $P_i$ and two halftwists with respect to $P_i$, then that fragment must be at lower level than the other fragments with endpoints on both ends of $P_i$ and at least one upper endpoint. 
        \item If $P_i$ is a lower cylinder, the fragments of $P_i$ with endpoints on only one end of $P_i$ and with at least one upper endpoint must be at lower level than the fragments they conflict with that have endpoints on both ends of $P_i$ and at least one upper endpoint in the spine sublist that contains the conflict. 
        \item If $P_i$ is a lower cylinder and has a fragment with upper endpoints on both ends of $P_i$ and two halftwists with respect to $P_i$ must be at lower level than the other fragments with endpoints on both ends of $P_i$ and at least one upper endpoint. 
        \item The previous four statements with ``lower level" replaced by ``higher level" are true if the fragments mentioned have at least one lower endpoint.
    \end{enumerate}
    
\end{lemma}

\begin{proof}

    In the first case, let an upper endpoint of the fragment $g$ with endpoints on both ends of $P_i$ lie in between the two upper endpoints of a fragment $f$ with endpoints on only one end of $P_i$. Then $f$ must be at lower level than $g$ by construction of $P_i$ being an upper cylinder.

    In the second case, if $P_i$ has a fragment $f'$ with two halftwists with respect to $P_i$ and an upper endpoint, then any other fragment attached to both ends of $P_i$ and having an upper endpoint has at most one halftwist with respect to $P_i$ and hence is at higher level than $f'$.

    Case 3 follows from case 1, and similarly case 4 follows from case 2 by applying a reflection with respect to the plane where the spine lies. Case 5 follows from cases 1 to 4 since if a fragment has a lower endpoint, then it has one more halftwist with respect to $P_i$ than a fragment with an upper endpoint.

\end{proof}

\begin{figure}
    \centering
    \includegraphics[width=0.5\textwidth]{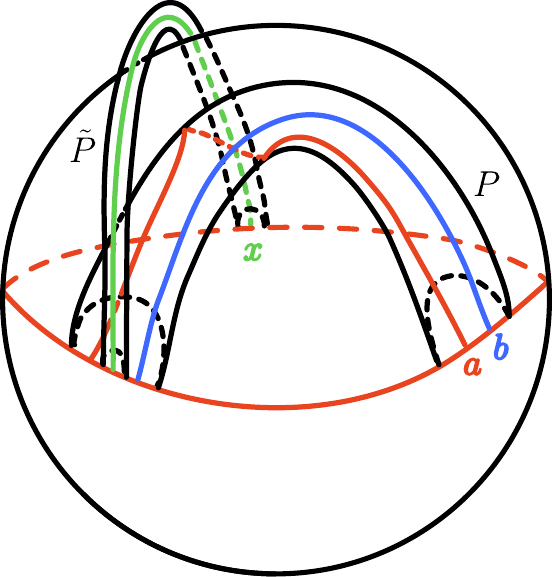}
    \caption{Illustration of the counterpositive statement of \cref{prop: in-between fragments lie above}.}
    \label{fig: in between fragments lie above}
\end{figure}

\begin{proposition}\label{prop: in-between fragments lie above}
   Suppose fragments $a$ and $b$, with $a$ at lower level than $b$ and both with at least one upper endpoint, are embedded in an upper cylinder $P$, with the disk sublists of $P$ restricted to $a$ and $b$ being $(a,b/a,b)$ or $(a,b,a/b)$. Suppose there exists a fragment $x$ such that an upper endpoint of $x$ lies on the spine between a pair of endpoints of $(a,b)$ in one disk sublist of $P$ (i.e. the disk sublists restricted to $a,b$ and such endpoint of $x$ are $(a,x,b/a,b)$, $(a,x,b,a/b)$ or $(a,b,x,a/b)$). If $x$ is at same or lower level than $a$, then $x$ is not embedded in an upper cylinder.
\end{proposition}
\begin{proof}
    We only prove the case where the disk sublists restricted to $a, b$ and the endpoint of $x$ lying on the spine between the pair of endpoints of $(a,b)$ as in the statement is of the form $(a,x,b/a,b)$. The other cases are similar. We prove the counterpositive: if $x$ is embedded in an upper cylinder, then $x$ is at higher level than $a$. (\Cref{fig: in between fragments lie above})
    
    If $x$ is embedded in $P$, then it must be at higher level than $a$, by \cref{lem: levels of fragments on cylinder}. 

    Consider the case where $x$ is embedded in another upper cylinder $\tilde{P}$. The upper endpoint of $x$ in between the pair of endpoints of $(a,b)$ lies on the intersection between the spine and the disk to which $P$ is glued to. Also, $P$ and $\tilde{P}$ are upper cylinders. Therefore $\tilde{P}$ crosses over $P$, which implies that $x$ crosses over $a$ or $b$.
\end{proof}

\begin{corollary}\label{cor: in-between fragments are above or below}
    Let $a$ and $b$ be fragments embedded in a cylinder $P$ with $a$ at lower level than $b$. Let $x$ be another fragment embedded in a cylinder $P'$ (not necessarily distinct from $P$) which has an endpoint lying on the spine between a pair of endpoints of $(a,b)$ in one disk sublist of $P$.
    \begin{enumerate}
        \item Assume $P$ is a lower cylinder and $a,b,x$ are lower endpoints. If $x$ is at same or higher level than $b$, then $x$ is not embedded in a lower cylinder.
        \item Assume $P$ is an upper cylinder and $a,b,x$ are lower endpoints. If $x$ is at same or higher level than $a$, then $x$ is not embedded in an upper cylinder. 
        \item Assume $P$ is a lower cylinder and $a,b,x$ are upper endpoints. If $x$ is at same or lower level than $b$, then $x$ is not embedded in a lower cylinder. 
    \end{enumerate}
\end{corollary}
\begin{figure}[h]
    \centering
    \includegraphics[width=0.5\textwidth]{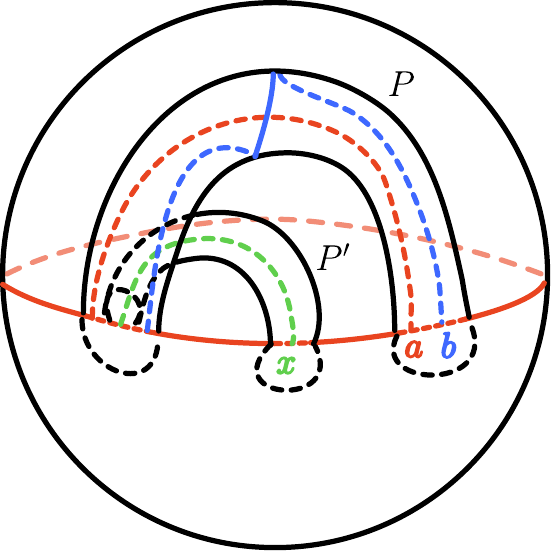}
    \caption{Illustration of the counterpositive statement of case 2 of  \cref{cor: in-between fragments are above or below}. Notice that one end of $P'$ is glued to $P$.}
    \label{fig: in-between fragments are above or below}
\end{figure}
\begin{proof}
    We only need to prove 2, because 1 follows from \cref{prop: in-between fragments lie above} and 3 follows from 2 by applying a reflection with respect to the plane where the spine lies. We prove the counterpositive: if $P$ and $P'$ are upper cylinders and the endpoints of $a,b,x$ are lower endpoints, then $x$ is at lower level than  $a$ (\Cref{fig: in-between fragments are above or below}). Let $P$ be an upper cylinder and $x$ have a lower endpoint in between two lower endpoints of $a$ and $b$ in the same disk sublist. If $x$ can also be embedded in $P$, then $x$ is at lower level than $b$ due to \cref{lem: levels of fragments on cylinder}.
    If $x$ is embedded in another upper cylinder $P'$, then it conflicts with either both or exactly one of $a$ and $b$. 
    In the latter case, the crossing number of $x$ with the fragment it does not conflict with must be 0, where the crossing number is the minimum number of crossing points between two fragments up to ambient isotopy. If $x$ conflicts with $a$, then $x$ is at lower level than $a$, because else $x$ would need to pierce the cylinder $P$ or cross over $P$, making it cross also fragment $b$; contradicting $x$ having crossing number 0 with $b$. 
    
    If $x$ conflicts only with $b$, it is not possible that $x$ is at higher level than $b$, because then $P'$ would cross over $P$, which implies that $x$ crosses over $a$ non-trivially. That is a contradiction of the fact that $x$ has crossing number 0 with $a$.
    
    The case where $x$ conflicts with both $a$ and $b$ implies that $x$ is at lower level than $a$. Indeed, it is not possible that $x$ is at a level in between the levels of $a$ and $b$, because $P'$ would need to pierce through $P$. If $x$ is at higher level than $b$, then $P'$ must cross over $P$. But that is not possible because $x$ has a lower endpoint in between two lower endpoints of $a$ and $b$. Therefore, $x$ needs first to cross under $a$ or $b$ before crossing over $b$, which is a contradiction to being leveled. Hence in this case $x$ is at lower level than $a$. 
\end{proof}

\subsection{Attaching cylinders to the surface}

\subsubsection{Conditions to avoid self-intersections of the surface}
In this section, given a hamiltonian leveled spatial graph $\G$, we investigate how to produce a surface $\mathcal{S}$ where $\G$ embeds such that $\mathcal{S}$ has no self-intersections. 

We examine when gluing a cylinder $P$ to $\mathcal{S}$ causes self-intersections. For that, we need to consider its interactions with already attached cylinders having endpoints in their spine sublist which, with respect to the spine list of $\G$, alternate with, are contained by, or contain endpoints in a spine sublist of $P$. Up to reflection of $\mathcal{S}$ with respect to the plane where the spine lies, four cases can arise at every endpoint: \begin{enumerate}
    \item upper cylinder with an upper endpoint and lower cylinder with a lower endpoint,
    \item upper cylinder with an upper endpoint and upper cylinder with a lower endpoint,
    \item upper cylinder with an upper endpoint and lower cylinder with an upper endpoint,
    \item upper cylinder with a lower endpoint and lower cylinder with an upper endpoint.
\end{enumerate}
We need to consider if the cylinders intersect each other only at the ends where their endpoints lie. If the cylinders do not intersect there, then it is possible to embed them without intersections.

We assume that we have a hamiltonian leveled spatial graph $\G$ consecutively embedded on a surface~$\mathcal{S}$.

\begin{definition}
    A cylinder $P$ and the endpoints of the fragments embedded in it \emph{interact} with another cylinder $P'$ and its endpoints if there exists an endpoint of a fragment of $P'$ which, in the spine list of $\G$, is in between two endpoints belonging to the same spine sublist of $P$ or vice versa. In particular, we say that $P$ \emph{encloses} $P'$ if a spine sublist of $P'$ is entirely in between two endpoints belonging to the same spine sublist of $P$. If a spine sublist of $P$ interacts with a spine sublist of $P'$ but neither encloses the other, then the endpoints of their fragments \emph{alternate}.
\end{definition}

\begin{figure}
    \centering
    \includegraphics[width=0.5\textwidth]{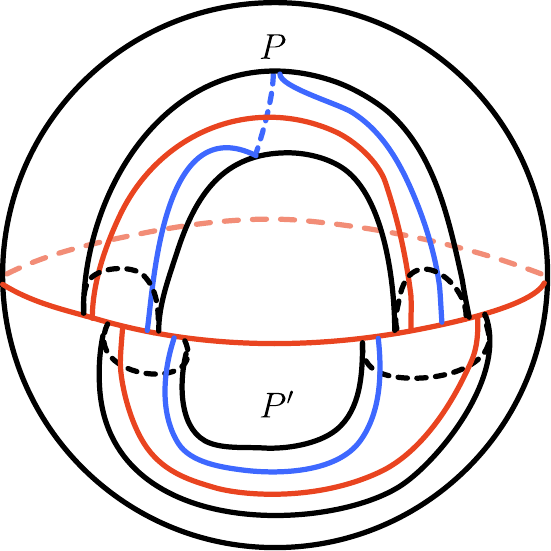}
    \caption{An upper cylinder $P$ with upper endpoints can freely interact with a lower cylinder $P'$ with lower endpoints.}
    \label{fig: interaction upper regular - lower regular}
\end{figure}

In case 1, the upper cylinder $P$ has an upper endpoint $e$ while the lower cylinder $P'$ has a lower endpoint $e'$. The ends of the two cylinders where $e$ and $e'$ lie do not intersect each other as the plane where the spine lies separates them. (\Cref{fig: interaction upper regular - lower regular}) 

We cover cases 2, 3 and 4 in \cref{prop: interaction upper irregular - upper regular}, \cref{prop: interaction lower irregular - upper regular} and \cref{cor: interaction lower irregular - upper irregular} respectively. 

\begin{figure}
    \centering
    \includegraphics[ width=0.5\textwidth]{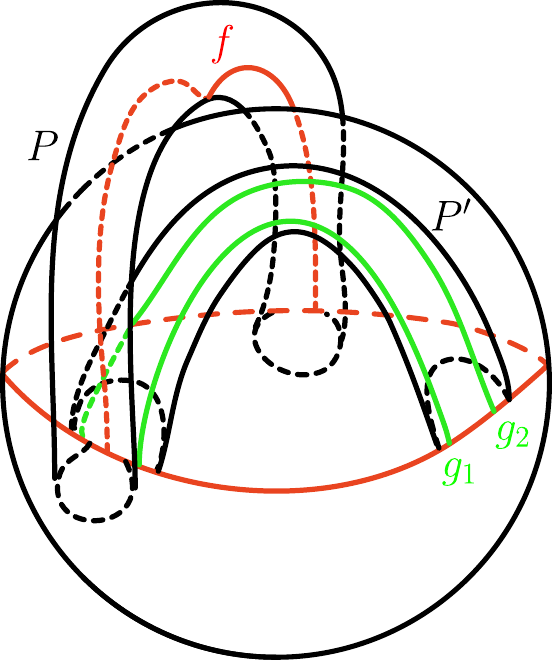}
    \caption{Illustration of the counterpositive statement of \cref{prop: interaction upper irregular - upper regular}.}
    \label{fig: prop about two upper cylinders with upper and lower endpoints}
\end{figure}

\begin{proposition}\label{prop: interaction upper irregular - upper regular}
    Let $P$ and $P'$ be upper cylinders. If the fragments $f_1, \dots, f_s$ of $P$ which conflict with a fragment of $g_1, \dots, g_t$ of $P'$ are at lower level than the fragments of $P'$ they conflict with, then the lower endpoints of $f_1, \dots, f_s$ do not interact with the upper endpoints of $g_1, \dots, g_t$ of $P'$ they conflict with. (See~\Cref{fig: prop about two upper cylinders with upper and lower endpoints} for the counterpositive.)
\end{proposition}
\begin{proof}
    We prove the case where $P'$ encloses $P$, arguing by contradiction; the other cases are similar.
    Suppose that the fragments $f_1, f_2$ of $P$ conflict with fragment $g_2$ of $P'$ and not with $g_1$. Moreover, suppose that the upper endpoints of $g_1$ and $g_2$ enclose the lower endpoints of $f_1, f_2$.

    $f_1$ and $f_2$ are at lower level than $g_2$. Therefore, if $P'$ encloses $P$, then $P$ needs to cross under $P'$, because $f_1, f_2$ have lower endpoints enclosed by upper endpoints of $g_1, g_2$. This implies $f_1, f_2$ cross under $g_1$, which contradicts $f_1$ and $f_2$ not conflicting with $g_1$.
\end{proof}

\begin{proposition}\label{prop: interaction lower irregular - upper regular}
    Let $P$ be a lower cylinder and $P'$ be an upper cylinder such that $P=P_j$ and $P'=P_i$ following the notation of \cref{def: consecutively embedded}. If $i>j$, then there does not exist a fragment $f$ of $P'$ which has an upper endpoint lying on the spine list in between two upper endpoints of fragments $g_1, g_2$ of $P$ belonging to the same spine sublist of $P$.
    If $i<j$, then $P$ and $P'$ do not intersect each other if and only if there exists a path from $g_1$ to $g_2$ along the boundary of the disk where $g_1,g_2$ lie which only crosses over $P'$ or below the fragments embedded in $\mathcal{S}_0$.
\end{proposition}
\begin{proof}
    We prove the counterpositive statement for the case $i>j$.
    
    By hypothesis the spine list restricted to the endpoints of the fragments $f, g_1, g_2$ on the interacting ends of $P$ and $P'$ is of the form $(g_1, f, g_2)$. If $i > j$, then the endpoint of $f$ in between those of $g_1$ and $g_2$ is contained in the boundary of $P$. Being $P$ a lower cylinder, this means that $f$ needs to cross under $g_1$ or $g_2$, which contradicts $P'$ being an upper cylinder.

    In case $i<j$, if there exists a path from the endpoint of $g_1$ to the endpoint of $g_2$ along the boundary of the disk with endpoints of $g_1, g_2$ on its boundary which only crosses over $P'$ or below the fragments embedded in $\mathcal{S}_0$, then it is possible to slightly move the path in order to form a disk $D$ to which the cylinder $P'$ can be glued. As the disk only crosses over $P'$ or below the fragments embedded in $\mathcal{S}_0$, then $P$ and $P'$ do not intersect.

    If instead there is no path from the endpoint of $g_1$ to the endpoint of $g_2$ along the boundary of the disk with endpoints of $g_1, g_2$ on its boundary which only crosses over $P'$ or below the fragments embedded in $\mathcal{S}_0$, meaning that every path from the endpoint of $g_1$ to the endpoint of $g_2$ crosses under $P'$ and above the fragments embedded in $\mathcal{S}_0$, then $P$ and $P'$ intersect by the intermediate value theorem applied to any path and fragment $f$.
\end{proof}

\begin{corollary}\label{cor: interaction lower irregular - upper irregular}
    Let $P$ be a lower cylinder and $P'$ be an upper one such that $P$ has an upper endpoint~$e$ and $P'$ has a lower endpoint~$e'$. Then the two cylinders can be placed on the resulting surface as long as they are not interacting.
\end{corollary}
\begin{proof}
    Follows from \cref{prop: interaction lower irregular - upper regular}. By symmetry, neither cylinder can enclose the other nor alternate with the other.
\end{proof}

\subsubsection{Conditions to ensure cellularity}
In this section, we investigate which hamiltonian leveled spatial graphs that are consecutively embedded are cellular. 
Recall that the subgraph of a consecutively embedded leveled spatial graph $\G$ embedded in a surface $S_i$, built successively by gluing the ends of cylinders to a sphere, is called $\G_i$.
\begin{lemma}\label{lem: embedding is cellular iff meridian and longitude are killed on cylinder}
    Let $\G_i$ be a hamiltonian leveled spatial graph consecutively embedded in a surface $\mathcal{S}_i.$ Suppose fragments $f_1, \dots, f_t$ of the spatial graph $\G_i$ are embedded in a cylinder $P_i$, while the spatial subgraph $\G_{i-1}\subset \G_i$ is cellular embedded in the surface $\mathcal{S}_{i-1}$ and $\G_i - \G_{i-1} = \{f_1, \dots, f_t\}.$
    \begin{enumerate}[label=(\arabic*)]
    \item If $P_i$ is glued to two different disks $D_1, D_2$ of the cellular embedding of $\mathcal{S}_{i-1}~\setminus~\G_{i-1}$, then the resulting embedding of $\G_i$ in $\mathcal{S}_i$ is cellular. 
    \item If $P_i$ is glued to a single disk $D$, then the resulting embedding is cellular if and only if there exists at least one conflicting pair of fragments in $f_1, \dots, f_t$ in the conflict graph with respect to the boundary of $D.$
    \end{enumerate}
\end{lemma}
\begin{proof}
    Assume that the fragment which is attached to both ends of $P_i$ is $f_j$.
    \begin{enumerate}
        \item In the case of $P_i$ being glued to two different disks, consider a closed path $\gamma_1 \in \mathcal{G}$ formed by $f_j$ and part of the spine joining the endpoints of $f_j$. All simple closed curves on $\mathcal{S}_i \setminus \G_{i-1}$ which do not bound a disk in $\mathcal{S}_i \setminus \G_{i-1}$ intersect $\gamma_1$ or the boundaries of $D_1$ or $D_2.$ This proves that gluing $P_i$ keeps the embedding cellular.
    
        \item Consider $P_i$ glued to $D$ as a torus $T$ with a marked point $t_0$, which coincides with the boundary of $D$. Two loops embedded in a torus form a cellular embedding if and only if their crossing number is greater than 0. In our case, loops are fragments and two fragments can cross each other only in $t_0$. Therefore, the embedding of $P_i$ in the torus $T$ with marked point $t_o$ is cellular if and only if there exist at least two fragments $f_1,f_2$ such that their endpoints alternate on the boundary of $D.$

    \end{enumerate}
\end{proof}

\begin{figure}
    \centering
    \includegraphics[width=0.8\textwidth]{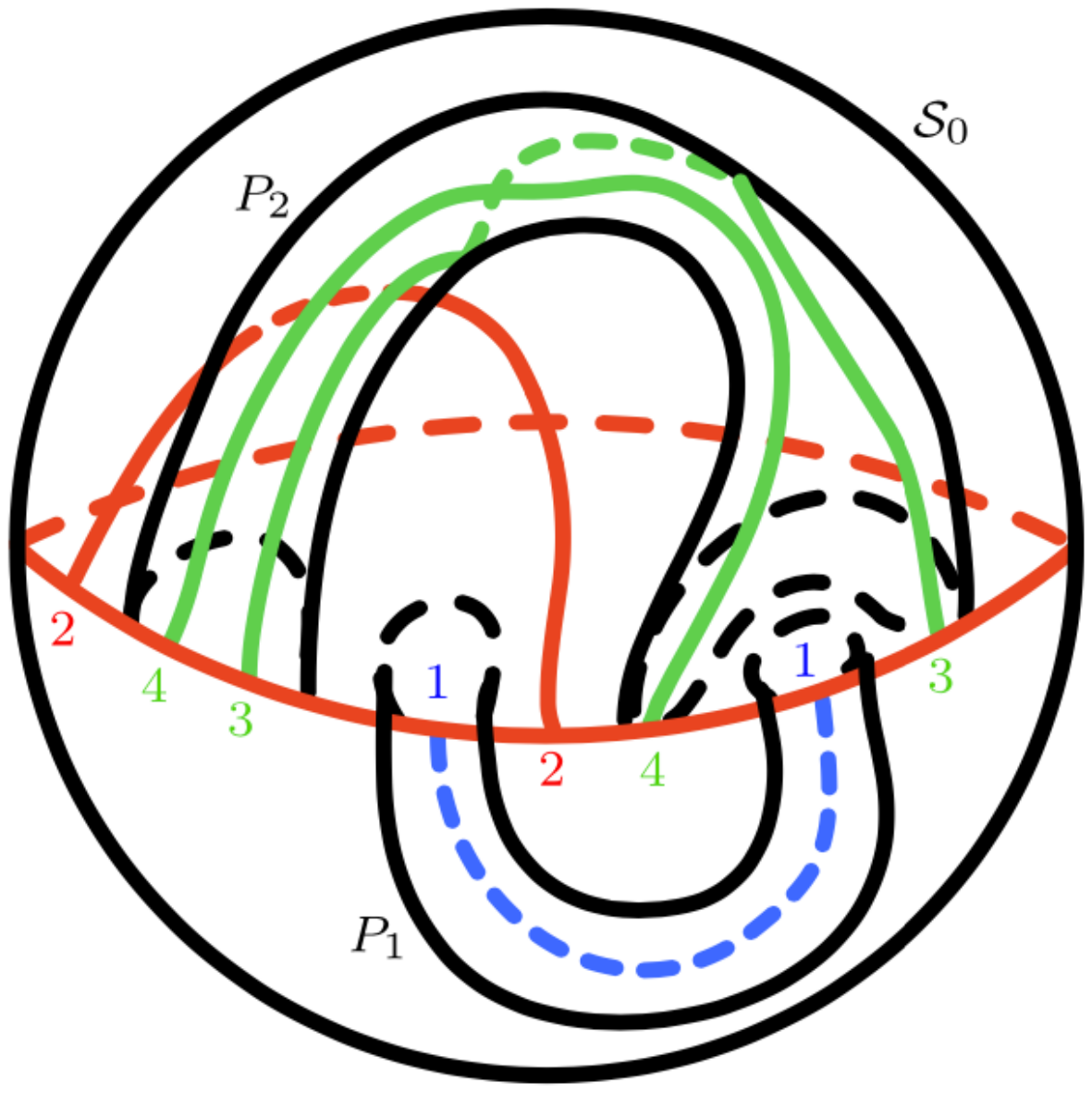}
    \caption{The embedding described in \cref{exm: conflict with respect to boundary of disk is important}.}
    \label{fig: conflict with respect to boundary of disk is important}
\end{figure}

\begin{example}\label{exm: conflict with respect to boundary of disk is important}
    Consider the leveled spatial graph $\G$ given by $(3, 1, 4, 2, 1, 3, 4, 2).$ We would like to embed $\G$ consecutively such that all endpoints are upper endpoints, fragment $2$ is embedded in the upper hemisphere of a sphere $\mathcal{S}_0$, fragment $1$ is embedded in a lower cylinder $P_1$ glued to the two disks formed by fragment $2$ (notice that this merges the two disks of the upper hemisphere into a single disk $D$) and fragments $3$ and $4$ are embedded together in an upper cylinder $P_2$ glued to $D$. See \Cref{fig: conflict with respect to boundary of disk is important} for a picture of the embedding.
    
    This embedding of $\G$ in $\mathcal{S}_2$ is not cellular, because fragments $3$ and $4$ do not conflict with respect to the boundary of $D$ by \cref{lem: embedding is cellular iff meridian and longitude are killed on cylinder}.
    Note in particular that $3$ and $4$ conflict with respect to the spine of $\G$, but this does not make the embedding cellular.
    
\end{example}

In \cref{prop: after gluing the surface has n-2+t disks} it is shown that the amount of disks obtained from gluing the cylinder $P_i$ to a surface $\mathcal{S}_{i-1}$ depends only on how many fragments are embedded in $P_i$, as long as the embedding is cellular.
\begin{proposition}\label{prop: after gluing the surface has n-2+t disks}
    Let $\G_i$ be a hamiltonian leveled spatial graph consecutively embedded in a surface $\mathcal{S}_i$. Let $\G_{i-1}$ be the (leveled) subgraph of $\G_i$ consecutively embedded in $\mathcal{S}_{i-1}$ and $\G_{P_i}$ be the subgraph of $\G_i$ embedded in $P_i$. Suppose that $\G_{i-1}$ is cellular embedded in $\mathcal{S}_{i-1}$ with $n$ disks. Suppose that $\G_{P_i}$ is cellular embedded in $P_i$ and $t$ is the number of fragments of $\G_{P_i} \subset \G_i$.
    
    Then, $\mathcal{S}_i \setminus \G_i$ has $n-2+t$ disks.
\end{proposition}
\begin{proof}
    The cellular embedding of $\G_{P_i}$ in $P_i$ has $t$ disks. Indeed, the first fragment of $\G_{P_i}$ which is attached to both ends of $P_i$ transforms an annulus into a disk, and every other fragment of $\G_{P_i}$ transforms an existing disk in two disks, adding one disk to the total per fragment.

    If $P_i$ is glued to two different disks of $\mathcal{S}_{i-1}$, then $\mathcal{S}_i \setminus \G_i$ has $n-2 +t$ disks. 
    
    If $P_i$ is glued to a single disk $D$ and the resulting embedding of $\G_i$ in $\mathcal{S}_i$ is cellular as seen in \cref{lem: embedding is cellular iff meridian and longitude are killed on cylinder}, then $\mathcal{S}_i \setminus \G_i$ has $(n-1) + (t-1)$ disks. Indeed, for each end of $P_i$, the disk of $P_i \setminus \G_{P_i}$ that contains $\partial P_i$ which is not part of the spine is merged with $D$. Notice that, due to cellularity, the disks merged on the two ends of $P_i$ are different.
\end{proof}

Given that the embedding is cellular, we adapt the ``Face tracing algorithm" of a cellular graph embedding explained in \cite{MR1855951} to be able to trace the boundaries of its disks.

\begin{fact}[Face tracing algorithm] \label{fact: face tracing algorithm}
     Let $\G$ be a leveled spatial graph with hamiltonian spine cellular consecutively embedded in a surface~$\mathcal{S}$. Consider the spine list of $\G$ and an endpoint~$e$ of a fragment~$f$. Go to the endpoint after $e$ on the spine list of $\G$ with same type as $e$ (i.e. upper if $e$ is an upper endpoint and lower otherwise). Go to the other endpoint of the fragment with endpoint~$e$ and repeat this process. The obtained closed walk is the boundary of a disk of the embedding. Starting at different endpoints gives all disks of the cellular embedding.
\end{fact}

\subsection{The algorithm for the surface construction} \label{sec: algorithm hamiltonian spine}
In this section, we provide an algorithm which, given as input a hamiltonian leveled spatial graph $\G$ in representative leveling by means of its spine list $\ell$, gives as output a surface $\mathcal{S}$ where $\G$ is embedded cellular and consecutively. The algorithm chooses two consecutive levels that are placed on a sphere and scans through all possible partitions of the remaining set of fragments of $\G$ into subsets, checking whether the fragments in each subset of a partition can be placed on a cylinder $P$ together. Partitions which do not yield a surface are identified using the results of the previous sections and are discarded.
The algorithm ends successfully if all subsets of one partition are placed on cylinders and ends unsuccessfully if all partitions are discarded.

The algorithm is extended to leveled spatial graphs with any spine, not necessarily hamiltonian, and to multi-leveled embeddings.

\begin{algorithm}\label{algo: hamiltonian leveled}

Let $\G$ be a hamiltonian leveled spatial graph with spine $C$, spine list~$\ell$ and $n$ levels. 

For $1\leq c < n $ choose the levels~$c$ and~$c+1$. Generate an embedded sphere $\mathcal{S}_0=S^2$ with the equator being the spine $C$ and  with fragments at level~$c$ being embedded in the lower hemisphere and fragments at level~$c+1$ in the upper hemisphere. Call this spatial subgraph~$\G_0$. 

Let $\mathcal{F}_{-\{c,c+1\}}$ be the set of $N$~fragments of $\G$ that are not in level~$c$ and~$c+1$. Generate all possible partitions $P_k, 1\leq k \leq \text{e}^{-1}\sum_{m=0}^{\infty}\frac{m^N}{m!}$ (Dobińsky's formula \cite{dobinsky1877}), of $\mathcal{F}_{-\{c,c+1\}}$ into disjoint nonempty subsets. That is, the elements of $P_k$ are subsets $P_k^j$, $1\leq j\leq J(k)$ of fragments, where $J(k)$ is the number of the subsets in $P_k$. The fragments in $P_k^j$ will be candidates for fragments embedded in one cylinder which by abuse of notation we also denote by $P_k^j$.

For $P_k$, there are $2^{J(k)}$ choices of labelling its elements with labels $\{U,L\}$.  
 This label indicates whether the cylinder $P_k^j$ is an upper or lower cylinder. (See \cref{exm: notation of partitions} for the notation.)

 The endpoints of the fragments in $P_k^j$ can be separated into two disk sublists. That is, for $P_k^j$ with $n$ endpoints, there are $\frac{1}{2}(n-1)n$ separations. Labelling the endpoints with $\{u,l\}$ results in $2^{n-1}(n-1)n$ labelled separations of $P_k^j$. This label indicates whether an endpoint is an upper or lower endpoint.

Enumerate the partition with labelled cylinders and labelled endpoints as $\mathcal{P}_i$. A labelled element in $\mathcal{P}_i$ is denoted by $P_i^j$. (See \cref{exm: assignments}.)

For $i$ from $1$ to the number of labelled partitions, consider $\mathcal{P}_i$.
For $0 \leq j \leq J(i)$, discard $\mathcal{P}_i$ if any of the following cases holds.

\begin{enumerate}
    \item $P_i^j$ has a disk sublist for which there exists no disk in $\mathcal{S}_{j-1} \setminus \G_{j-1}$ whose boundary contains all endpoints in the disk sublist. 

\item The list of endpoints in $P_i^j$  is not one of those in \cref{prop: fragments on cylinder}.

\item $P_i^j$ is an upper (resp. lower) cylinder, its fragments with endpoints on only one end of $P_i^j$ and with at least one upper endpoint are at higher level than the fragments of $P_i^j$ they conflict with that have endpoints on both ends of $P_i^j$ and at least one upper endpoint. (\cref{lem: levels of fragments on cylinder})

\item $P_i^j$ is an upper (resp. lower) cylinder with a fragment with endpoints on both ends of $P_i^j$, two halftwists with respect to $P_i^j$ and at least one upper endpoint which is at higher level than the other fragments with endpoints on both ends of $P_i^j$ and at least one upper endpoint. (\cref{lem: levels of fragments on cylinder})

\item $P_i^j$ is an upper (resp. lower) cylinder, its fragments with endpoints on only one end of $P_i^j$ and with at least one lower endpoint are at lower level than the fragments of $P_i^j$ they conflict with that have endpoints on both ends of $P_i^j$ and at least one lower endpoint. (\cref{lem: levels of fragments on cylinder})

\item $P_i^j$ is an upper (resp. lower) cylinder with a fragment with endpoints on both ends of $P_i^j$, two halftwists with respect to $P_i^j$ and at least one lower endpoint which is at lower level than the other fragments with endpoints on both ends of $P_i^j$ and at least one lower endpoint. (\cref{lem: levels of fragments on cylinder})

\item There exists a cylinder $P_i^k$ in $\mathcal{P}_i$ that has the same label as $P_i^j$. And for a pair of fragments $F_s$ at level $s$ and $F_t$ at level $t>s$ in $P_i^j$ that are conflicting with respect to the spine, there exists in $\ell$ an upper (resp. lower) endpoint $x$ at level lower (resp. higher) than $s$ of a fragment $g\in P_i^k$ that lies between two upper (resp. lower) endpoints of $F_s$ and/or $F_t$ that belong to the same spine sublist of $P_i^j$. (\cref{prop: in-between fragments lie above}, \cref{cor: in-between fragments are above or below})

\item Tracing the disks which the fragments of $P_i^j$ have endpoints on using \cref{fact: face tracing algorithm} reveals that the endpoints of all fragments of $P_i^j$ lie on a single disk $D$ and no pair of fragments of $P_i^j$ conflicts with respect to the boundary of $D$. (\cref{lem: embedding is cellular iff meridian and longitude are killed on cylinder})

\item There exists a cylinder $P_i^h$ with $h<j$ in $\mathcal{P}_i$, such that one of the following cases holds:
\begin{itemize}
    \item Both $P_i^h$ and $P_i^j$ are upper (resp. lower) cylinders, and some lower (resp. upper) endpoints of $P_i^h$ interact with some upper (resp. lower) endpoints of $P_i^j$, and the fragments of $P_i^h$ interacting are at lower (resp. upper) level than the interacting fragments of $P_i^j$ they conflict with. (\cref{prop: interaction upper irregular - upper regular})
    \item $P_i^h$ is a lower (resp. upper) cylinder, and $P_i^j$ is an upper (resp. lower) one, and there exists a fragment of $P_i^j$ which has an upper (resp. lower) endpoint in between two upper (resp. lower) endpoints of fragments of $P_i^h$ belonging to the same spine sublist of $P_i^h$. (\cref{prop: interaction lower irregular - upper regular})
    \item $P_i^h$ is an upper (resp. lower) cylinder, and $P_i^j$ is a lower (resp. upper) one, and there exists a fragment of $P_i^h$ which has an upper (resp. lower) endpoint in between two upper (resp. lower) endpoints of fragments $g_1, g_2$ of $P_i^j$ belonging to the same spine sublist of $P_i^j$ such that there is no path from an endpoint of $g_1$ to an endpoint of $g_2$ along the boundary of the disk with endpoints of $g_1, g_2$ on its boundary which only crosses over (resp. under) $P_i^h$ or below (resp. above) the fragments embedded in $\mathcal{S}_0$.~\footnote{To detect if a path from an endpoint of $g_1$ to an endpoint of $g_2$ along the boundary of the disk $D$ with endpoints of $g_1, g_2$ on its boundary only crosses over a cylinder $P_i^h$, it is enough to trace the boundary of $D$ with \cref{fact: face tracing algorithm} and check if there exists a path from the endpoint of $g_1$ to the endpoint of $g_2$ which includes at least one fragment embedded on a cylinder $P_i^k$ with $k > h$ and no fragments embedded on $\mathcal{S}_l$ for $l <i$ which conflict with $f$.} (\cref{prop: interaction lower irregular - upper regular})
    \item $P_i^h$ is an upper (resp. lower) cylinder with a fragment with an upper (resp. lower) endpoint interacting with a fragment with a lower (resp. upper) endpoint of $P_i^j$ and $P_i^j$ is a lower (resp. upper) cylinder. (\cref{cor: interaction lower irregular - upper irregular})
    \end{itemize}
\end{enumerate}

If $\mathcal{P}_i$ has not been discarded, embed the fragments in $P_i^j$ in a cylinder $P_i^j$ that is upper or lower according to the label given to $P_i^j$. Glue $P_i^j$ to the disk(s) of $\mathcal{S}_{j-1} \setminus \G_{j-1}$ which were traced in case~4, obtaining the surface $\mathcal{S}_j.$ Call $\G_j$ the subgraph of $\G$ induced by $\G_{j-1}$ and the fragments in $P_i^j$. Consider the set $P_i^{j+1} \in \mathcal{P}_i$ and repeat the process.

If $P_i^{j+1}$ does not exist, then the algorithm ends successfully.

If $\mathcal{P}_i$ is discarded, consider $\mathcal{P}_{i+1}$. If $\mathcal{P}_{i+1}$ does not exist, set $c=c+1$. If $c=n-1$, the algorithm ends unsuccessfully.

\end{algorithm}

\begin{example}\label{exm: notation of partitions}
Consider three fragments $A, B, C$, that is $N=3$. Then there are \mbox{$\text{e}^{-1}\sum_{m=0}^{\infty}\frac{m^3}{m!}= 5$} partitions $P_1 = \{\{A\},\{B\},\{C\}\}$, $P_2 = \{\{A\},\{B,C\}\}$, $P_3 = \{\{B\},\{A,C\}\}$, $P_4 = \{\{C\},\{A,B\}\}$, \mbox{$P_5 = \{\{A,B,C\}\}$}. That is, $J(1)=3$, $J(2)=J(3)=J(4)=2$, $J(5)=1$.\\
Some examples of cylinders are $P_1^3 = \{C\}$, $P_2^1 = \{A\}$, $P_2^2 = \{B,C\}$, $P_3^1 = \{B\}$. \\
There are $2^{J(2)}= 2^2=4$ possibilities to label the cylinders of $P_2$ as upper or lower cylinders: $\{\{A\}^U,\{B,C\}^U\}$, $\{\{A\}^U,\{B,C\}^L\}$, $\{\{A\}^L,\{B,C\}^U\}$, $\{\{A\}^L,\{B,C\}^L\}$.
\end{example}

\begin{example}\label{exm: assignments}
If a cylinder contains a single edge $P_i^j=\{A\}$ on level~$b$ the spine sublist of the fragments in $P_i^j$ is $(A_b^1, A_b^1)$. Since $P_i^j$ contains only two endpoints, there is only one possibility of dividing the spine sublist into disk sublists, namely $(A_b^1/ A_b^1)$. That leaves eight combinations of labelling the endpoints as well as the cylinder, of which four are labeled as upper cylinders $P_i^{j,U}$ and four as lower cylinders $P_i^{j,L}$.

That is, if the leveled graph has only three levels and there is only one fragment on level~$b=3$, for $c=1$ the algorithm steps through the following partitions:  \begin{align*} \mathcal{P}_1&=\{(A_{3,u}^{1} /A_{3,u}^{1})^U\},& \mathcal{P}_2&=\{(A_{3,u}^{1} / A_{3,l}^{1})^U\},& \mathcal{P}_3&=\{(A_{3,l}^{1}/A_{3,u}^{1})^U\}, \\ \mathcal{P}_4&=\{ (A_{3,l}^{1}/A_{3,l}^{1})^U\}, &
 \mathcal{P}_5&=\{(A_{3,u}^{1}/A_{3,u}^{1})^L\},&
\mathcal{P}_6&=\{(A_{3,u}^{1}/A_{3,l}^{1})^L\}, \\
\mathcal{P}_7&=\{(A_{3,l}^{1}/A_{3,u}^{1})^L\},& \mathcal{P}_8&=\{ (A_{3,l}^{1}/A_{3,l}^{1})^L\} . \end{align*}
\end{example}

The following theorem is a consequence of the results of the previous sections.
\begin{theorem}\label{thm: algorithm hamiltonian leveled}
    Let $\G$ be a hamiltonian leveled spatial graph. If \cref{algo: hamiltonian leveled} applied to $\G$ ends successfully, its output is an oriented closed surface $\mathcal{S} \subset \mathbb{R}^3$ without self-intersections where $\G$ is cellular embedded.
\end{theorem}

Algorithm~\ref{algo: hamiltonian leveled} requires the spine to be a hamiltonian cycle. However, we can generalize the algorithm to general (non-hamiltonian) leveled embeddings.

If we assume that the spine of the leveled spatial graph is a non-hamiltonian cycle, fragments do not need to be edges, but they all are planar connected graphs by definition. In order to place handles under such spatial fragments, we need to attach $n$-punctured spheres instead of $2$-punctured spheres (i.e. cylinders) to build the surface. Note that an $n$-punctured sphere can be obtained by gluing $n-1$ cylinders together.

\begin{algorithm}\label{algo: non-hamiltonian algorithm}
Let $\G$ be a leveled spatial graph with spine~$C$. Let $e$ be any edge of $C$ and choose a spanning tree~$T$ of $\G$ which contains $C- e$. The choice of $T- (C- e)$ is independent of the choice of $e$.

Embed $T \cup e$ in a sphere~$\mathcal{S}_0$, such that $C$ is the equator of $\mathcal{S}_0$. This embedding is cellular, because $T \cup e$ is a connected trivial spatial graph. Contracting all edges in $T - C$ gives a hamiltonian leveled spatial graph $\G'$ to which \cref{algo: hamiltonian leveled} can be applied. If the algorithm terminates successfully, a cellular embedding $(\G', \mathcal{S}')$ is obtained. The edge contraction can be reversed by vertex splitting and subdividing edges in $\mathcal{S}'$ to obtain a cellular embedding $(\G, \mathcal{S})$.

\end{algorithm}

The following theorem is a consequence of \cref{thm: algorithm hamiltonian leveled} and \cref{algo: non-hamiltonian algorithm}.

\begin{theorem}\label{thm: extended algorithm leveled}
    Let $\G$ be a leveled spatial graph. If \cref{algo: non-hamiltonian algorithm} applied to $\G$ ends successfully, its output is an oriented closed surface $\mathcal{S} \subset \mathbb{R}^3$ without self-intersections where $\G$ is cellular embedded.
\end{theorem}

\begin{figure}[h!]
    \centering
    \includegraphics[width=0.8\textwidth]{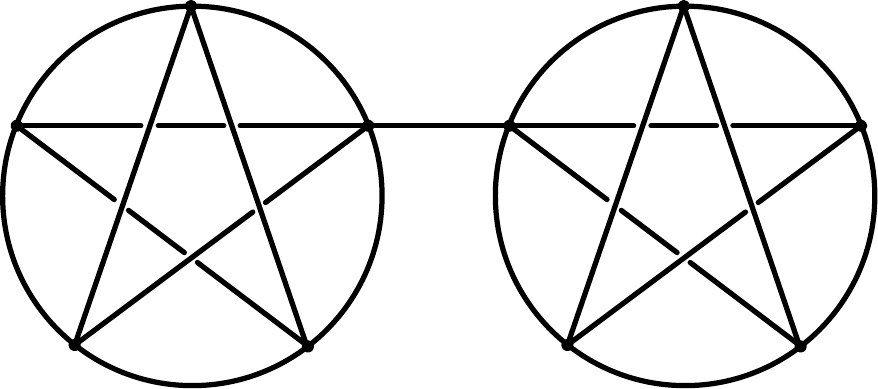}
    \caption{A multi-leveled embedding of the graph $K_5 \_ K_5$ with two blocks, each isomorphic to~$K_5$, and one planar graph, isomorphic to $K_2$.}
    \label{fig: multileveled embedding}
\end{figure}

Note that the underlying abstract graph of a leveled embedding consists of at most one non-planar 2-connected component, by definition. For example, the graph $K_5 \_ K_5$, obtained by joining two copies of $K_5$ by an edge (see \Cref{fig: multileveled embedding}), does not admit a leveled embedding, because it contains two non-planar 2-connected components.
However, if multiple leveled embeddings are attached to each other via trivial spatial graphs, the algorithms presented before can still be applied. We call this generalization of leveled embeddings \emph{multi-leveled embeddings}:

\begin{definition} A connected spatial graph $\G$ that can be constructed as described in the following is \textbf{multi-leveled}.\\
    Consider a tree~$T$ and the bipartite partition $X_1,X_2$ of the vertices of $T$. To every vertex $v_1^i \in X_1$ assign a leveled spatial graph~$\G_i$, $1\leq i \leq |X_1|$. The $\G_i$ are called \textbf{blocks} of the embedding. To every vertex $v_2^j \in X_2$, $1\leq j \leq |X_2|$, assign a connected planar graph~$H_j$ with at least $\deg(v_2^j)$~vertices. For each vertex~$v_1^i$ adjacent to $v_2^j$, identify some vertices of $H_j$ with some vertices of the block~$\G_i$ such that identified vertices belong to exactly one planar graph and one block. (See \Cref{fig: example complicated multi-leveled embedding} and \Cref{fig: tree multi-leveled embedding})

    \begin{figure}[h!]
        \centering
       \includegraphics[width=0.9\textwidth]{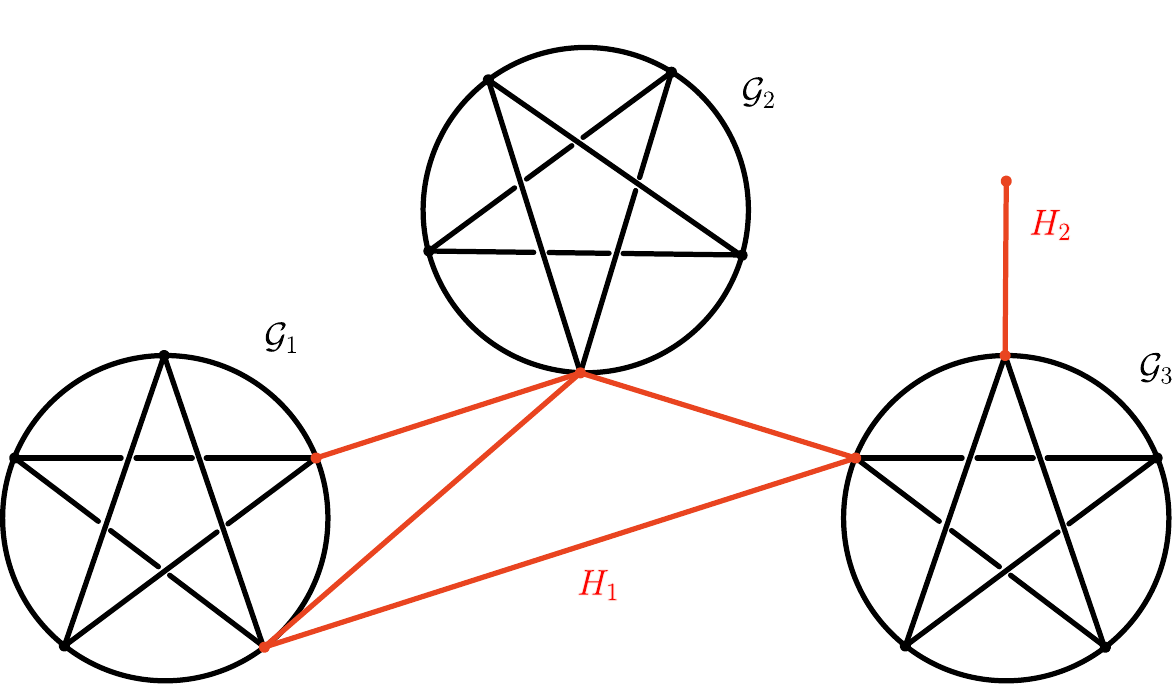}
        \caption{A multi-leveled embedding with three blocks and two connected planar graphs.}
        \label{fig: example complicated multi-leveled embedding}
    \end{figure}
    \begin{figure}[h!]
        \centering
        \includegraphics[width=0.7\linewidth]{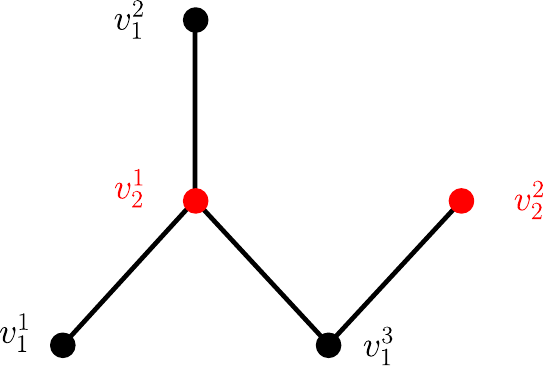}
        \caption{The tree associated to the embedding of \Cref{fig: example complicated multi-leveled embedding}.}
        \label{fig: tree multi-leveled embedding}
    \end{figure}
 
\end{definition}

\begin{theorem}
    Let $\G$ be a spatial graph with a multi-leveled embedding with blocks $\G_1, \dots, \G_n$. If: \begin{enumerate}[label=(\arabic*)]
        \item \cref{algo: non-hamiltonian algorithm} applied to $\G_i$ produces a surface $\mathcal{S}_i$ for each $1 \le i\le n$,
        \item every $H_j$ can be cellular embedded in a $d_j$-punctured sphere~$S_{d_j}$ such that the vertices of $H_j$ identified with vertices $w_j^1, \dots, w_j^t$ of block $\G_i$ all lie in the boundary of a puncture of $S_{d_j}$,
        \item $w_j^1, \dots, w_j^t$ all lie on the boundary of a single disk $D_{i,j}$ of $\mathcal{S}_i \setminus \G_i$,
        \item the disks can be chosen such that $D_{i,j}\neq D_{i',j'}$ for $i,j \neq i', j'$,
      
    \end{enumerate} then there exists a closed oriented surface $\mathcal{S}$ such that $(\G, \mathcal{S})$ is cellular.
\end{theorem}
\begin{proof}
    The surface $\mathcal{S}$ is constructed by attaching the punctured sphere $S_{d_j}$ to $\mathcal{S}_i$. Therefore the boundary of the puncture where vertices of $H_j$ are identified with vertices of $\G_i$ and the boundary of the cell $D_{i,j}$ are identified and the interior of $D_{i,j}$ removed. The embeddings of $H_j$ in $S_{d_j}$ and of $\G_i$ on $\mathcal{S}_i$ are cellular by assumption. Moreover, condition $(3)$ together with \cref{lem: embedding is cellular iff meridian and longitude are killed on cylinder} ensure that the embedding $(\G, \mathcal{S})$ is also cellular.
\end{proof}

\subsubsection{Results of the algorithm}

\begin{figure}[h!]
    \centering
    \includegraphics[width=0.7\textwidth]{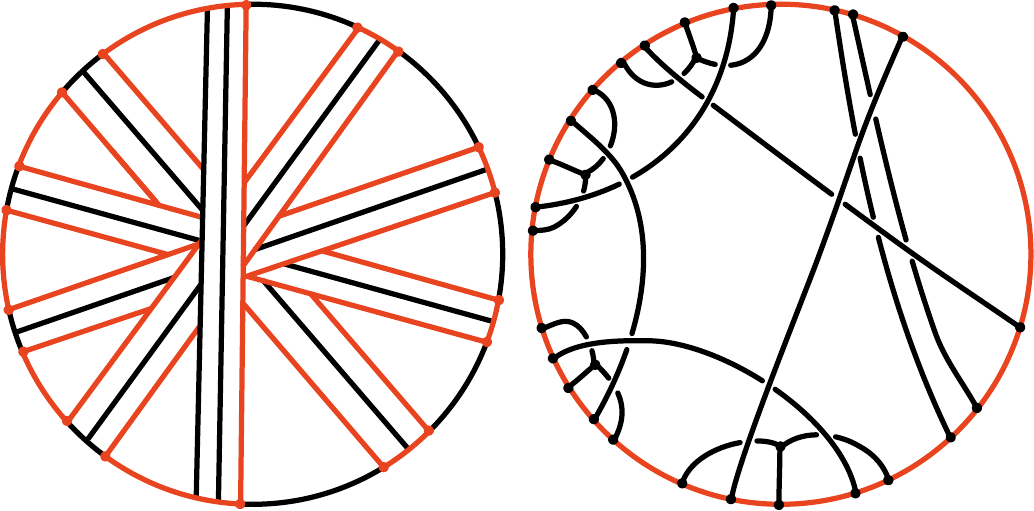}
    \caption{Left: the new spine chosen for $\G$ in red. Right: the leveled embedding with the new spine.}
    \label{fig: changing the spine it can be cellular embedded}
\end{figure}

\begin{figure}[h!]
    \centering
    \includegraphics[width=0.7\textwidth]{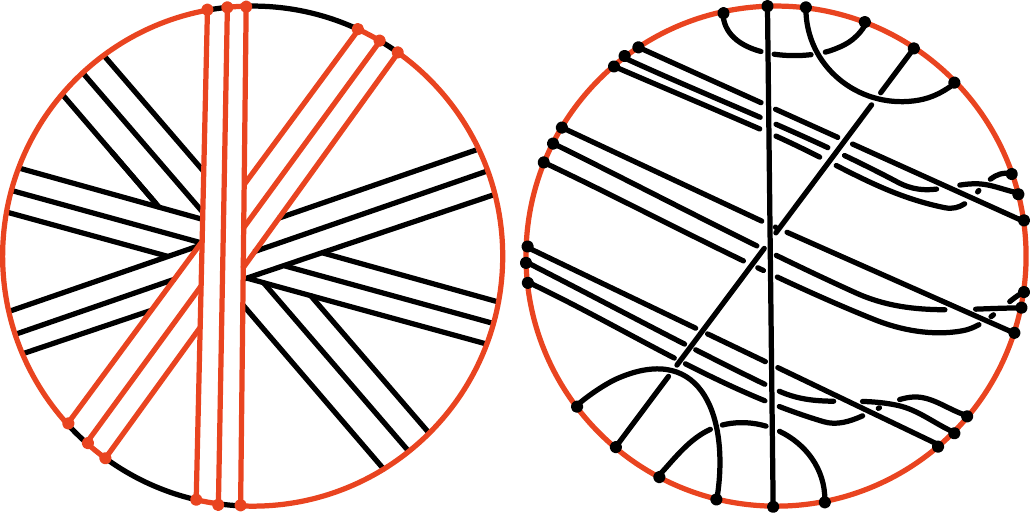}
    \caption{On the left in red, the new hamiltonian spine chosen for $\G$. On the right, the leveled embedding with the new spine.}
    \label{fig: changing the spine to hamiltonian it can still be cellular embedded}
\end{figure}

The algorithm gives a sufficient condition for a leveled spatial graph to be cellular embedded in an oriented surface. To investigate in how far the successful termination of the algorithm is also a necessary condition, we tried to construct a leveled spatial graph $\mathcal{G}$ for which the algorithm ends unsuccessfully for all possible choices of the spine. We have not been able to construct such an example. In all graphs we studied, it was always possible to find a spine for which \cref{thm: algorithm hamiltonian leveled}, \cref{thm: extended algorithm leveled}, \cref{prop: few levels can always be cell embedded}, or \cref{cor: up to 5 levels may be cell embedded} apply. 
As an example, consider the spatial graph $\G$ described in \cref{exm: graph not embeddable with the algorithm unless you change spine}. Although for the given spine the algorithm does not find a surface $\mathcal{S}$ where $\G$ cellular embeds, it is possible to choose a non-hamiltonian spine for $\G$, shown in \Cref{fig: changing the spine it can be cellular embedded}, for which \cref{algo: non-hamiltonian algorithm} terminates successfully. Furthermore, a different hamiltonian spine for $\G$ can be chosen, shown in \Cref{fig: changing the spine to hamiltonian it can still be cellular embedded}, which yields a leveled embedding with five levels. Corollary~\ref{cor: up to 5 levels may be cell embedded} guarantees the existence of a surface where $\G$ cellular embeds, and indeed \cref{algo: hamiltonian leveled} terminates successfully if applied to this leveled embedding.
Note that it follows from this example that all hamiltonian leveled spatial graphs with a spine list of the form $(1^1, 1^2, 1^3, 2^1, 2^2, 2^3, \dots, n^1, n^2, n^3, 1^3, 1^2, 1^1, 2^3, 2^2, 2^1, \dots, n^3, n^2, n^1)$, for $n \in \mathbb{N}$, can be cellular embedded by applying \cref{algo: hamiltonian leveled} to a different choice of the spine, analogously to the example.

\begin{figure}[h!]
    \centering
    \includegraphics[width=0.7\textwidth]{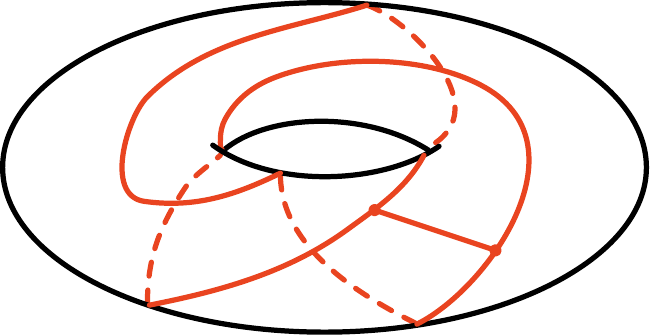}
    \caption{A trefoil knot with unknotting tunnel cellular embeds in a torus but does not admit a leveled embedding.}
    \label{fig: trefoil cellular embeds on torus}
\end{figure}

It is important to highlight that admitting a leveled embedding for a spatial graph $\G$ is not a necessary condition for the existence of a surface $\mathcal{S}$ where $\G$ cellular embeds. Indeed, the trefoil knot with an unknotting tunnel is an example of a spatial graph that does not admit a leveled embedding, because for every unknotted cycle chosen as spine, the fragment will cross the spine non-trivially. It does though have a cellular embedding in a torus, as shown in \Cref{fig: trefoil cellular embeds on torus}.

In \cite{3PageEmbeddingsContinua}, it is shown that any spatial graph $\G$ can be embedded in a 3-book, where the vertices lie on the back of the 3-book and edges may be contained in several leaves (i.e. the pages of the 3-book). Placing a vertex on every edge of $\G$ which intersects the back of the 3-book and adding an edge between two consecutive vertices on the back which are not already adjacent yields a leveled embedding of a spatial graph $\G'$ in three levels such that $\G \subset \G'$. In other words, every spatial graph has a super-graph with a leveled embedding for which \cref{algo: non-hamiltonian algorithm} ends successfully.
 

\section{Outlook}

For a given spatial graph $\G$, if \cref{algo: non-hamiltonian algorithm} ends unsuccessfully for every choice of spine $C$ of a leveled embedding of $\G$, it does not imply that $\G$ cannot be cellular embedded in any surface. However, in every example we studied, it is possible to choose a different spine for which \cref{algo: non-hamiltonian algorithm} ends successfully. This motivates the following conjectures.
\begin{conjecture}\label{conj: algorithm is necessary}
    If a spatial graph $\G$ has a leveled embedding and admits a cellular embedding in some surface, there is a choice of spine for which \cref{algo: non-hamiltonian algorithm} ends successfully.
\end{conjecture}
\begin{conjecture}\label{conj: every leveled can be cell emb}
    If a spatial graph $\G$ has a leveled embedding, then a cycle can be found that is the spine of a leveled embedding of $\G$ for which \cref{algo: non-hamiltonian algorithm} ends successfully.
\end{conjecture}

A similar conjecture can be stated for hamiltonian leveled spatial graphs. The situation shown in \Cref{fig: changing the spine to hamiltonian it can still be cellular embedded} is representative for all examples of hamiltonian leveled spatial graphs that we studied.
\begin{conjecture}
    If a spatial graph $\G$ has a hamiltonian leveled embedding, then there exists a choice of hamiltonian spine and corresponding hamiltonian leveled embedding of $\G$ for which \cref{algo: hamiltonian leveled} ends successfully.
\end{conjecture}

The algorithms~\ref{algo: hamiltonian leveled} and~\ref{algo: non-hamiltonian algorithm}, while giving a constructive procedure to find a surface $\mathcal{S}$ where a spatial graphs cellular embeds, are intractable to actually implement for spatial graphs with a high number of levels. For a practical implementation, some heuristics would be needed to optimize the algorithm, based for example on which choice of cylinders is more likely for the algorithm to end successfully. 

Cellular embeddings of graphs are related to circular embeddings. A circular embedding is a cellular embedding $(\G, \mathcal{S})$ where the boundaries of the faces of $\mathcal{S}\setminus \G$ are cycles of $\G$. A famous open conjecture in graph theory states that every 2-connected graph has a circular embedding~\cite{JAEGER19851}. In case of a 2-connected hamiltonian leveled spatial graph $\G$ for which \cref{algo: hamiltonian leveled} ends successfully, it is possible to determine whether the resulting cellular embedding is circular by applying \cref{fact: face tracing algorithm} to all faces of the cellular embedding: if every edge of $\G$ appears in the boundary of exactly two faces, then the embedding is circular. The embedding $(\G,\mathcal{S})$ is ensured to be circular if no cylinder attachments merge disks in the construction of $\mathcal{S}$ as performed in \cref{algo: hamiltonian leveled}.

Investigating which embeddings obtained using \cref{algo: non-hamiltonian algorithm} are circular would give insights in the circular embedding conjecture for 2-connected graphs with a multi-leveled embedding.

\section{Acknowledgement}
The authors want to thank Joanna A. Ellis-Monaghan for discussions that motivated this research.

\section{Bibliography}

\bibliographystyle{alpha}
\bibliography{bibliography}

\end{document}